\documentclass[twoside]{article}
\usepackage{eurosym}
\usepackage[utf8]{inputenc}
\usepackage{graphicx}
\usepackage{latexsym}
\usepackage{amsmath}
\usepackage{amsfonts}
\usepackage{amssymb}
\usepackage{hyperref}
\usepackage[usenames]{color}

\newcommand{\disp}{\displaystyle}
\newcommand{\cA}{{\mathcal A}}
\newcommand{\rA}{{\rm A}}
\newcommand{\rB}{{\rm B}}
\newcommand{\bB}{{\bf B}}
\newcommand{\cB}{{\mathcal B}}
\newcommand{\cC}{{\mathcal C}}
\newcommand{\rC}{{\rm C}}
\newcommand{\rD}{{\rm D}}
\newcommand{\cD}{{\mathcal D}}
\newcommand{\rE}{{\rm E}}
\newcommand{\EE}{\mathbb{E}}
\newcommand{\be}{{\bf e}}
\newcommand{\bF}{\mathbb{F}}
\newcommand{\rF}{{\rm F}}
\newcommand{\cF}{{\mathcal F}}
\newcommand{\rG}{{\rm G}}
\newcommand{\cG}{{\mathcal G}}
\newcommand{\HH}{\mathbb{H}}
\newcommand{\rH}{{\rm H}}
\newcommand{\bH}{{\bf H}}
\newcommand{\rI}{{\rm I}}
\newcommand{\cJ}{{\mathcal J}}
\newcommand{\rK}{{\rm K}}
\newcommand{\cL}{{\mathcal L}}
\newcommand{\rL}{{\rm L}}
\newcommand{\cM}{{\mathcal M}}
\newcommand{\rM}{{\rm M}}
\newcommand{\cO}{{\mathcal O}}
\newcommand{\PP}{{\mathbb{P}}}
\newcommand{\cP}{{\mathcal P}}
\newcommand{\rP}{{\rm P}}
\newcommand{\cQ}{{\mathcal Q}}
\newcommand{\rQ}{{\rm Q}}
\newcommand{\RR}{{\mathbb{R}}}
\newcommand{\R}{\RR}
\newcommand{\rR}{{\rm R}}
\newcommand{\rS}{{\rm S}}
\newcommand{\rT}{{\rm T}}
\newcommand{\rV}{{\rm V}}
\newcommand{\cV}{{\mathcal V}}
\newcommand{\rW}{{\rm W}}
\newcommand{\rX}{{\rm X}}
\newcommand{\fin}{\hfill\mbox{$\quad{}_{\Box}$}}
\newcommand{\fineq}{\vspace{-.75cm$\fin$}\par\bigskip}
\newcommand{\fineqnum}{\vspace{-.4cm$\fin$}\par\bigskip}
\newcommand{\pe}[2]{\langle #1,#2\rangle}
\newcommand{\n}[1]{\disp \|#1\|}
\newtheorem{theo}{\bf \sffamily Theorem}
\newtheorem{prop}{\bf \sffamily Proposition}
\newtheorem{rem}{\bf \sffamily Remark}
\newtheorem{lemma}{\bf \sffamily Lemma}
\newtheorem{definition}{\bf \sffamily Definition}

\begin{document}

\title{{\Large \bfseries\sffamily} Stochastic energy balance
	climate models with Legendre weighted diffusion and  a cylindrical Wiener process forcing}
\author{\bfseries\sffamily G. D\'{\i}az \& J.I. D\'{\i}az \thanks{Partially supported the UCM Research Group MOMAT (ref. 910480) and  the projects MTM2017-85449-P and PID2020-112517GB-I00 of the DGISPI, Spain. \hfil\break \indent {\sc
			Keywords}: Stochastic Energy Balance Model, climate diffusion, cylindrical Wiener process, parabolic Legendre diffusion, stochastic attractor, random dynamical systems. 
		}}
\date{}
\maketitle

\centerline{\em Dedicated to Georg Hetzer on occasion of his 75th birthday}

\begin{abstract}
We consider a class of one-dimensional nonlinear stochastic parabolic problems associated to Sellers and Budyko diffusive energy balance climate models with a Legendre weighted diffusion and an additive cylindrical Wiener processes forcing. Our results use in an important way that, under suitable assumptions on the Wiener processes, a suitable change of variables leads the problem to a pathwise random PDE, hence an essentially "deterministic" formulation depending on a random parameter. Two applications are also given: the stability of solutions when the Wiener process converges to zero and  the asymptotic behaviour of solutions for large time.
\end{abstract}

\section{Introduction}
We consider a class of stochastic of climate diffusive
energy balance models of the type 
$$
(\rE_{\beta ,\rG})\left\{ 
\begin{array}{l}
\dfrac{\partial u}{\partial t}-\dfrac{\partial }{\partial x}\left (\big (
1-x^{2})\dfrac{\partial u}{\partial x}\right  )+g(u)\in 
\rQ\rS(x)\beta (u)+f+\rG\dfrac{\partial 
	\rW}{\partial t},\\[0.2cm] 
u(x,0)=u_{0}(x),%
\end{array}
\right. 
$$
where $x\in \rI\doteq (-1,1)$ and $t>0$.  
Notice that, in fact, $u=u(x,t;\omega)$ and $\dfrac{\partial 
\rW}{\partial t}=~\dfrac{\partial  \rW}{\partial t}(x,t;\omega)$, where $x\in \rI,~t\geq 0$ and $\omega$ is in the probability space $\{\Omega, \cF,\PP\}$. In which follows we will use the notation $\varrho (x)=1-x^{2},~ x\in \rI$. We mainly assume $\rQ>0$,
\par
(\bH$_{g}$) $g:\RR\rightarrow \RR$ is a continuous increasing
function.
\par
(\bH$_{\beta }$) $\beta $ is a {\em bounded} maximal monotone graph in $%
\RR^{2},~$i.e.$~\beta (s)\in [m,\rM]$, $\forall s\in 
\RR$.
\par
(\bH$_{s}$) $\rS:\rI\rightarrow \RR,~\rS\in 
\rL^{\infty }(\rI)$, $\rS_{1}\geq \rS(x)\geq 
\rS_{0}>0$ $\ \mbox{a.e.}~x\in \rI$.

\noindent On the data $u_{0}$ and $f$ we will assume different conditions which will be presented later. The expression $\rG\dfrac{\partial \rW}{\partial t}$ denotes a time dependent {\em white or real noise} which we always assume to be corresponding to a {\em cylindrical Wiener processes}, as we will explain later.
\par
This kind of problems where proposed by R. North and R.F. Cahalan in 1982 (\cite{NCa}) for the modeling of non-deterministic variability (as, for
instance, the volcanoes actions) in the context of diffusive energy balance climate
models. We recall that the distribution of temperature $u(x,t)$ is expressed
pointwise after some averaging processes. The spatial variable $x$
is given by $x={\rm sin}~\theta $ with $\theta $ the latitude on a
supposed spherically symmetric Earth, which leads to the above mentioned Legendre weighted diffusion operator. The deterministic equation on a
general representation of the Earth as a Riemannian manifold without
boundary was carried out in \cite{DT} and \cite{Diaz-Hetzer}, among other references.
\par
Very often, the natural degenerate diffusion given by the Legendre weighted diffusion operator $\dfrac{\partial }{\partial x}\left (\big (1-x^{2}\big )\dfrac{\partial u}{\partial x}\right )$ is sometimes replaced by the usual $1d$-Laplacian operator and, then,  the absence of boundary
conditions for the degenerate diffusion arising in $(\rE_{\beta ,\rG})$ is replaced by asking homogeneous Neumann boundary conditions (since in $(\rE_{\beta ,\rG})$ the meridional heat flux $(1-x^{2})\dfrac{\partial u}{\partial x}$ vanishes at the poles $x=\pm 1$). In that paper we will not
use such a simplification, improving some previous results in the
literature. From the modeling point of view, the balance of energies leads
to the problem 
$$
\left\{ 
\begin{array}{ll}
\dfrac{\partial u}{\partial t}-\dfrac{\partial }{\partial x}\bigg (\varrho 
\dfrac{\partial u}{\partial x}\bigg )=\rR_{a}-\rR_{e}, & \quad
\hbox{in $\rI\times (0,+\infty )$}, \\[0.2cm] 
u(x,0)=u_{0}(x), & \quad x\in \rI,
\end{array}
\right. 
$$
where the terms $\rR_{a}$ and $\rR_{e}$ must be specified by
means of constitutive laws (see, $e.g$., the monographs and surveys \cite{NK}, \cite{Diaz-Hetzer-Tello}, \cite{Kapper-Engler}, \cite{Hetzer Houston},  and \cite{Imkeller}). The {\em absorbed energy} $\rR_{a}$ depends, in a fundamental way, on the planetary {\em co-albedo} $\beta $ representing the
fraction of the incoming radiation flux which is absorbed by the surface. In
ice-covered zones, reflection is greater than over oceans, therefore, the
co-albedo is smaller. So, there is a sharp transition between zones of high
and low co-albedo. In the energy balance climate models, a main change of the
co-albedo occurs in a neighborhood of a critical temperature for which ice
become white, usually taken as $u=-10^{\circ }\rC$. In the so called 
{\em Budyko model} the different value of the co-albedo is modeled by means
of a discontinuous function of the temperature (\cite{Bu}). As usual in
PDEs, this function can be understood in the more general context of the
maximal monotone graphs in $\RR^{2}$. In particular, we assume that 
\begin{equation}
\beta (u)=\left\{ 
\begin{array}{ll}
\beta _{i}, & \quad \hbox{if $u<-10$,} \\[0.1cm] 
\left[ \beta _{i},\beta _{w}\right] , & \quad \hbox{if $u=-10$}, \\[0.1cm] 
\beta _{w}, & \quad \hbox{if $u>-10$},
\end{array}
\right.  \label{eq:BetaHeaviside}
\end{equation}%
where $m=\beta _{i}$ and $\rM=\beta _{w}$ represent the co-albedo in
the ice-covered zone and the free-ice zone, respectively and $0<\beta
_{i}<\beta _{w}<1$ (the value of these constants has been estimated by
observation from satellites). In contrast to the above assumption, in the so
called {\em Sellers model} (\cite{Sellers}) $\beta $ is assumed to be a
more regular function (at least, Lipschitz continuous) piecewise linear
function far from a neighborhood of $u=-10$. In both models, the whole
absorbed energy is given by $\rR_{a}(x,t,u)=\rQ\rS(x)\beta \big (u(x,t)\big )$ where $\rS(x)$ is the {\em insolation function} and $\rQ$ is the so-called {\em solar constant}.
\par
The Earth's surface and atmosphere, warmed by the Sun, emits part of the
absorbed solar flux as an infrared long-wave radiation. This energy $\rR_{e}$ is given by the Stefan-Boltzman law (when temperature is given in
Kelvin degrees) $\rR_{e}(x,t,u)=\mu (x,t)u(x,t)^{4}$, for
some $\mu (x,t)>\varepsilon _{0}>0,$ but, following the proposal by Budyko,
sometimes it is enough to consider a linearization of this law leading to
expressions of the type $\rR_{e}(x,t,u)=\mathrm{B} u(x,t)-f(x,t)$.
Obviously, $\rR_{e}$ includes the action due to the {\em greenhouse effect}. But any representation of $\rR_{e}$ is incomplete without
taking into account weather fluctuations on the thermal field which does not
obey to deterministic laws such as, for instance, the influence of gases
resulting from the eruption of volcanoes and their influence on the
emissivity coefficient. At the present time, according the Global Volcanism
Program of the Smithsonian Institution (GVP-SI), there are more than 1.500
active volcanoes and a very high number of dormant and extinct volcanoes. In
terms of probabilities, it is like thinking of a die with more than 1500
faces. The correct modeling can be formulated by adding some white or real noise as a external forcing.
Here, following previous studies (see, $e.g.$ \cite{NCa}) we will
assume that the stochastic term is given by some cylindrical Wiener process 
\begin{equation}
\rR_{e}(x,t,u)=g\big (u(x,t)\big )-f(x,t)+\rG(x,t)\dfrac{
	\partial \rW}{\partial t}(x,t),  \label{Def R-e}
\end{equation}%
since, essentially the stochastic variability of volcanoes affects only to
the time variable and not to the space variable. For a different study dealing with the impact of applying stochastic forcing to the Sellers energy balance climate model in the form of a fluctuating solar irradiance see \cite{Lucarini}.
\par
As said before, the absorbed energy is assumed to be of the form 
$$
\rR_{a}(x,t,u)\doteq \rQ\rS(x)\beta \big (u(x,t)%
\big ), 
$$
under assumptions (\bH$_{\beta }$) and (\bH$_{s}$). We recall the important
difference in the assumptions made on $\beta (u)$ in the so called Budyko
and Sellers models.
\par
We point out that some of the results of this paper remain valid for a
larger class of non necessarily cylindrical Wiener processes, as presented
in \cite{DD}. Nevertheless, in this paper it is assumed the above mentioned
simplifying version of the diffusion operator and moreover the constitutive
expression for $\rR_{e}$ following the suggestion made in \cite{Bu}
for the linear choice $g(r)=\rC r$ with $\rC$ a positive
constant. In the previous paper \cite{Diaz Langa Valero} problem $(\rE_{\beta ,\rG})$ was consider under additional conditions (besides the
above mentioned simplifications on the diffusion and the function $g(r)$: it
was assumed that $\rG(x,t)\rW(x,t)\doteq\rG(x)\rB(t)$ with $\rB(t)$ a Brownian motion).
\par
The main goal of this paper is to generalize the previous contributions by
the authors, \cite{DD} and~\cite{Diaz Langa Valero}, to the case in which
problem $(\rE_{\beta ,\rG})$ is formulated in a unifying
suitable class of cylindrical Wiener process. We also will provide of some
new results (see, $e.g.$ estimate \eqref{eq:comparisonestimate} bellow). We
point out that some different approach with some numerical experiences can
be found in the climatological studies \cite{NCa} and \cite{NK}.
\par
We recall that in the deterministic case ($\rG=0$) the
existence of solutions was given in \cite{Diaz Escorial} (see~\cite{DT} for
the generalization to bidimensional models) and that in the case in which $%
\beta $ is multivalued it was shown that there is lack of uniqueness of
solutions except in the class of the so called, {\em non degenerate solutions}. As it was shown in \cite{DD}, a curious fact is produced for
problem $(\rE_{\beta ,\rG})$: the presence of a stochastic
perturbation produces a kind of uniqueness of the solutions result
associated to any given monotone (univalued and discontinuous) section $b$
of the maximal monotone graph $\beta $. A similar comment emphasizing the
presence of the space-time white noise was quoted by \cite{GP2} in a general
context.
\par
The common point of view in our treatment of the stochastic problem 
$(\rE_{\beta ,\rG})$ is to consider the problem as a special
case of a class of abstract Cauchy problems, with an additive white noise,
on the Hilbert space $\rH$. The keystone of our treatment is based on a
basic idea, which under some different formulations, seems to be already
quoted in \cite{BenTem}  (see also \cite{DaPZ}, \cite{Arnold}, \cite{Caraballo-Langa-Valero 2003}, \cite{Diaz Langa Valero}, \cite{Bar} and \cite{BarRooper}).It consists in introducing a suitable change of variables reducing the Cauchy problems to a pathwise random PDE, hence an essentially "deterministic" formulation
depending on a random parameter $\omega \in \Omega $. Some authors call to this technique to reduce a Stochastic Differential Equation (respectively a Stochastic Partial Differential Equation) to a random Ordinary Differential Equation (respectively a random Partial Differential Equation) as the {\em Doss-Sussmann method} (introduced, independently, on 1977, in \cite{Doss} and \cite{Suss}). Nevertheless,  the special change we will introduce is not exactly the same (for instance in \cite{HK} this change of variables is connected with the Ornstein-Uhlenbeck equation, which is not our case)  and, in fact, was already used in \cite{BenTem} on 1972. 
\par
In order to present the assumptions on the stochastic part it is convenient to introduce some useful notation. We start by recalling (\cite{Diaz Escorial}) that for the treatment of the diffusion operator it is convenient to work with the {\em energy space} given by the weighted space
$$
\rV=\big \{w\in\rL^{2}(\rI):~w'\in \rL^{2}(\rI;\varrho)\big \}
$$
where 
$$
\rL^{2}(\rI:\varrho)=\left \{v:~\int _{\rI}\varrho |v|^{2}dx<+\infty\right \},$$
equipped with its norm
$$
\n{v}_{\rL^{2}(\rI:\varrho)}=\left (\int _{\rI}\varrho |v|^{2}dx\right )^{\frac{1}{2}}
$$
(we recall that $\varrho(x)=1-x^{2}$).
The treatment is made on the Hilbert pivot space $\rH=\rL^{2}(\rI)$ equipped with its usual norm 
$$
\n{u}_{\rL^{2}(\rI)}=\left (\int_{\rI}|u|^{2}dx \right)^{\frac{1}{2}}.$$
Notice that $\rV$ is a separable Hilbert space related to the norm
$$
\n{w}_{\rV}=\n{w}_{\rL^{2}(\rI)}+\n{w'}_{\rL^{2}(\rI:\varrho)}.
$$
Then we treat the diffusion operator by means of the functional operator $\cA:~\rV\rightarrow \rV'$ given by 
\begin{equation}
\cA u\doteq -\dfrac{\partial }{\partial x}\left( \varrho 
\dfrac{\partial }{\partial x}u\right),\quad u\in\rV.
\label{eq:leadingpart}
\end{equation}
It is also useful to define its realization as an  operator on  $\rH$, $\rA:~\rD(\rA)\rightarrow \rH$ with
\begin{equation}
\left \{
\begin{array}{l}
\rD(\rA)=\big \{v\in \rH:~\cA v\in \rH\big\},\\ [.2cm]
\rA v=\cA v,\quad \hbox{if $v\in \rD(\rA)$.}
\end{array}
\right .
\label{eq:leadingpartHoperator}
\end{equation}
On the other hand, it is well known (see, $e.g.$,  \cite{NK}), that $\rH=\rL^2(\rI)$  admits a Hilbertian base given by the eigenvectors $\be_{n}\in\rD(\rA)$ of the operator $\rA$ defined through the  orthonormal Legendre polynomials of degree $n$
\begin{equation}
\left \{
\begin{array}{l}
\be_{n}(x)=\sqrt{\dfrac{2}{2n+1}}\rP_{n}(x),\quad -1\le x\le 1,\\ [.4cm]
\mu_{n}=n(n+1)\quad\hbox{({\em the corresponding eigenvalues}).}
\end{array}
\right .
\label{eq:baseA}
\end{equation}
The polynomials are solutions of the Legendre equation and are given by
$$
\rP_{n}(x)=\dfrac{1}{2^{n}n!}\dfrac{d^{n}}{dx ^{n}}\big (x^{2}-1\big )^{n}. 
$$
Thus $\rP_{0}(x)\equiv 1,~\rP_{1}(x)=x,~\rP_{2}(x)=\dfrac{3x^{2}-1}{2},~\rP_{3}(x)=\dfrac{x(5x^{2}-3)}{2},\ldots $ and so on (see \cite{NSU}).
\par
In this, and other previous papers in the literature on this kind of problems, the considered  stochastic processes $\rG(x,t)\dfrac{\partial \rW}{\partial t}(x,t)$ can be classified in two different types:
\par
\noindent {\bf a)} the so-called {\em finite dimensional stochastic noise} (considered by most of the authors interested in asymptotic behaviour) given by
\begin{equation}
\rG(x,t)\dfrac{\partial \rW(x,t)}{\partial t}\doteq \pe{\Phi (x,t)}{\dfrac{d \bB (t)}{dt}},
\quad (x,t)\in\rI\times (0,\infty)
\label{eq:finitestochasticnoise}
\end{equation}
where $\Phi (\cdot,t)\in\big (\rV)^{m}$ satisfying sometimes some additional regularity conditions and where $\bB(t)=\{\rB^{1}(t),\ldots ,\rB^{m}(t)\}$ is a set of real mutually independent standard Brownian motions on a filtered probability space $\{\Omega ,\cF, \{\cF_{t}\}_{t\ge 0}, \PP \}$ (see Section 4 below), 
\par
\noindent {\bf b)} the so-called {\em infinite dimensional stochastic noise} (considered by some  authors interested in general existence and uniqueness of solutions results), in which now
\begin{equation}
\rG(x,t)\dfrac{\partial \rW(x,t)}{\partial t}\doteq \sum_{n\ge 1} \dfrac{1}{\sqrt{\mu_{n}}}\rG(x,t)\dfrac{d \rB^{n}(t)}{dt}\be_{n}(x)
\quad (x,t)\in\rI\times (0,\infty).
\label{eq:infinitestochasticnoise}
\end{equation}
This corresponds to assume
\begin{equation}
\rW(x,t)\doteq\sum_{n\ge  1}\dfrac{1}{\sqrt{\mu_{n}}}\rB^{n}(t)\be_{n}(x),\quad (x,t)\in\rI\times (0,+\infty),
\label{eq:Hypo Wiener cylindrical}
\end{equation}%
(see \eqref{eq:baseA}) where $\{\rB^{n}(t)\}_{n\ge 1}$ are mutually
independent Brownian motions on a filtered probability space $\{\Omega ,\cF, \{\cF_{t}\}_{t\ge 0},
\PP \}$, with right-continuous filtration.  Notice that a simple version of the, so called, {\em Wiener isometry} holds:
\begin{equation}
\EE\big [\n{\rW(\cdot,t)}^{2}_{\rH}\big ]=t\sum_{n\ge 1}\dfrac{1}{\mu_{n}}<+\infty,
\label{eq:simpleiosmetry}
\end{equation}
(see \cite{LiRo} and the Appendix). Then, the correct treatment of the term $\rG(.,t)$ is formulated as an operator
\begin{equation}
\rG(\cdot,t):~\rI\rightarrow \cL_{2}(\rH),\quad 0<t<\rT<\infty,
\label{eq:infinitestochasticnoiseG}
\end{equation}
where $\cL_{2}(\rH)$ denotes the Hilbert-Schmidt space on $\rH$ (see \cite[Appendix B]{LiRo} and the Appendix below for some details).
In that case, we can generate a function (denoted in the same way)
\begin{equation}
\rG(x,\cdot):~(0,\rT)\times \Omega\rightarrow \RR,\quad x\in\rI,
\label{eq:Gpredictable}
\end{equation}
such that it satisfies the condition which defines the so called {\em predictable processes}: it is $\cB([0,\rT])\otimes \cF$-measurable. 
\par

In the rest of the paper we will follow the usual notation in stochastic processes concerning the t-dependence (for instance $u(x,t)\doteq u_{t}(x)$). Connecting our problem with a more abstract setting, we are in the conditions in which we start with a given Gelfand triple 
$$
\rV\hookrightarrow \rH\hookrightarrow \rV^{\prime }
$$
and assume given an operator $\cA:\rV\rightarrow \rV^{\prime }$. Then, the problem $(\rE
_{\beta ,\rG})$ becomes a special case of the abstract stochastic Cauchy problem 
\begin{equation}
\left\{ 
\begin{array}{l}
du_{t}+\cA u_{t}dt=(\rR_{a}-\rR_{e})_{t}dt+\rG_{t}d\rW_{t},\quad t>0, \\%
[0.1cm]
u_{t}\big |_{t=0}=u_{0}\in \rH.
\end{array}
\right.   \label{eq:SDEabs}
\end{equation}
\par
Different kind of notions of solutions are possible and then it is crucial to formulate correctly the different assumptions on the data. At least formally, the problem  \eqref{eq:SDEabs} is equivalent to the integral identity 
\begin{equation}
u_{t}=u_{0}-\int_{0}^{t}(\cA u_{s}-f_{s})ds+\int_{0}^{t}\rG_{s}d\rW_{s},\quad \forall t>0,
\label{eq:strongsolution} 
\end{equation}
and $\{u_{t}\}_{t\geq 0}$  must be, at least, an adapted random process to the filtration satisfying in some sense the integral representation \eqref{eq:strongsolution}. We note that in a strong way, it would require that $\cA u_{t}-f_{t},~t\in (0,\rT)$, let $\PP$-$a.s$ Bochner integrable (see \cite[Appendix~A]{LiRo}). 
\par
On the other hand, for the finite dimensional stochastic noises we also assume in this paper \eqref{eq:finitestochasticnoise} and
\begin{equation}	
\int_{0}^{\rT}\n{\Phi_{s}} _{\rL_{2}(\rI)}^{2}ds<+\infty .	
\label{eq:paraisometriafinite}
\end{equation}
Then we introduce the notation 
\begin{equation}
(\Phi\cdot\rB)_{t}\doteq \int_{0}^{t}\pe{\Phi_{s}}{d\bB_{s}}
\label{def Zedafinite},
\end{equation}
whose {\em Wiener isometry} is
\begin{equation}
\EE\big [\n{(\Phi\cdot\rB)_{t}}^{2}_{\rL^{2}(\rI)}\big ]=\int_{0}^{t}\n{\Phi _{s}} _{\rL_{2}(\rI)}^{2}ds<+\infty.
\label{eq:Wienerisometryintegralfinite}	
\end{equation}
Notice that, in this special case and under the above conditions, we have $(\mathrm{G}\cdot \mathrm{W})\in C([\mathrm{0,T}]:\mathrm{V})$.
For the infinite dimensional stochastic noises \eqref{eq:infinitestochasticnoise} we assume \eqref{eq:infinitestochasticnoiseG},\eqref{eq:Gpredictable} and
\begin{equation}	
\int_{0}^{\rT}\n{\rG_{t}\cJ} _{\cL_{2}(\rH)}^{2}dt<+\infty ,	
\label{eq:paraisometria}
\end{equation}
for which the stochastic integral is well defined taking values in $\rH$, where $\cJ:\rH\rightarrow \rH$ is given by the property $\cJ\be_{n}=\dfrac{1}{\sqrt{\mu_{n}}}\be_{n}$. Then it is useful to introduce the notation
\begin{equation}
(\rG\cdot\rW)_{t}\doteq \int_{0}^{t}\rG_{s}d\rW_{s} 
\label{def Zeda}.
\end{equation}
Now the {\em Wiener isometry} is
\begin{equation}
\EE\big [\n{(\rG\cdot\rW)_{t}}^{2}_{\rL^{2}(\rI)}\big ]=\int_{0}^{t}\n{\rG_{s}\cJ} _{\cL_{2}(\rH)}^{2}ds<+\infty.
\label{eq:Wienerisometryintegral}	
\end{equation}
\par
We always will assume in this paper the additional condition
\begin{equation}
(\rG\cdot\rW)_{t}\in \cC([0,\rT];\rH)\cap \rL^{2}(0,\rT;\rV).  
\label{eq:Hypo W}
\end{equation}
We notice that when $\rG$ is  time independent then $(\rG\cdot\rW)_{t}=\rG\rW_{t}$.

We present now three notions of solutions which are relevant for the rest of the paper. By simplicity, we will follow the notations appearing in  \cite[Appendix G]{LiRo}. The reader can find in this reference an exposition about the mutual connections between the different notions of solutions.
\begin{definition}  A $\rV$-valued $\{u_{t}\}_{t\ge 0}$ predictable process, (as in \eqref{eq:Gpredictable}),  is called an ``analytically strong solution''  if the identity \eqref{eq:strongsolution} holds $\PP-a.s$ for each $t\in[0,\rT]$.
\end{definition}
In this framework, concerning the perturbation term $(\rR_{a}-\rR_{e})$, we assume
\begin{equation}
\rR_{a}-\rR_{e}\in \rL^{2}(0,\rT;\rV'). 
\label{eq:Hypo f(t,x)}
\end{equation}
Therefore, the identity  \eqref{eq:strongsolution} must be understood to be taking place on $\rV'$, but then the term $\cA u_{t}-(\rR_{a}-\rR_{e})_{t}$ must be in $\rV$, which is very restrictive in practice. In order to avoid this difficulty it is convenient to weaken the notion of solution. So, by taken duality products we get 
\begin{equation}
\pe{u_{t}}{v}=\pe{u_{0}}{v}+\int_{0}^{t}\pe{\sqrt{\varrho }\frac{
		\partial {u}_{s}}{\partial x}}{\sqrt{\varrho }\frac{\partial {v}}{\partial x}}ds+\int_{0}^{t}\pe{ (\rR_{a}-	\rR_{e})_{s}}{v}ds+\int_{0}^{t}\pe{\rG_{s}^{*}v}{d
	\rW_{s}}_{\rV'\times \rV}
\label{eq:AnalWeakSol}
\end{equation}
for $t>0$ and any $v\in\rV$.

Notice that for $v\in\rH$ one has the identity
$$
\int_{0}^{t}\pe{\rG_{s}v}{d\rW_{s}}_{\rV'\times \rV}=\int_{0}^{t}\pe{\rG_{s}^{*}v}{d\rW_{s}}_{\rH},
$$
where the right hand side term is well defined. In the case of the infinite dimensional stochastic noises the map $\widetilde{\rG}_{t}\doteq\pe{\rG_{s}^{*}v}{w}_{\rH}$ belongs $\cL_{2}(\rH;\RR)$ and
$$
\n{\widetilde{\rG}_{t}}_{\cL_{2}(\rH;\RR)}^{2}=\sum_{n\ge 1}|\pe{\rG_{t}^{*}v}{\be_{n}}_{\rH}|^{2}=\n{\rG_{t}^{*}v}_{\rH}^{2}<+\infty.
$$
Moreover, the relative Wiener isometry becomes
$$
\EE\bigg [\left |\int_{0}^{t}\pe{\rG_{s}^{*}v}{d\rW_{s}}_{\rH}\right |^{2}\bigg ]=\int^{t}_{0}
\n{\cJ^{*}\rG_{s}^{*}v}_{\rH}^{2}ds <\infty.
$$
Then we can introduce the following notion of solution:
\begin{definition}  A $\rH$-valued predictable process $\{u_{t}\}_{t\ge 0}$ is called an ``analytically weak solution'' if the identity  \eqref{eq:AnalWeakSol} holds for each $t\in[0,\rT]$ and  $v\in \rV$.
\end{definition}

Finally, a third notion of solution can be introduced when we assume 
\begin{equation}
\rR_{a}-\rR_{e}\in \rL^{2}(0,\rT;\rH) 
\label{eq:HypofH}
\end{equation}
and consider the operator $\rA:\rD(\rA)\rightarrow \rH$ defined in \eqref{eq:leadingpartHoperator}. We know  that $\rA$ generates a semi-group $\{\rS(t)\}_{t\ge 0}$ of contractions on $\rH$  (see \cite{Diaz Escorial}). Then we can introduce the following notion: 
\begin{definition}  A $\rH$-valued predictable process $\{u_{t}\}_{t\ge 0}$ is called a "mild solution" if
\begin{equation}
u_{t}=\rS(t)u_{0}+\int^{t}_{0}\rS(t-s)(\rR_{a}-\rR_{e})_{s}ds+\int^{t}_{0}\rS(t-s)\rG_{s}d\rW_{s},
\label{eq:mildsolution}
\end{equation}
holds for each $t\in[0,\rT]$.
\end{definition}
The concrete formulation of the change of variable we propose in this paper
consists in to consider a new unknown process defined on $\{\Omega ,\cF,\{\cF_{t}\}_{t\ge 0},\PP\}$ by 
\begin{equation}
y_{t}=u_{t}-(\rG\cdot\rW)_{t}.  \label{eq:Change of variables}
\end{equation}%
Then the stochastic Cauchy problem \eqref{eq:SDEabs} becomes the 
\textit{deterministic} Cauchy problem depending on the  random parameter, $\omega \in \Omega$,
\begin{equation}
\left\{ 
\begin{array}{l}
\dfrac{\partial y_{t}}{\partial t}-\dfrac{\partial }{\partial x}\big (
\varrho \dfrac{\partial y_{t}}{\partial x}\big )+\rR_{a}(\cdot
,t,y_{t}+(\rG\cdot \rW)_{t})=\rR_{e}(\cdot ,t,y_{t} +\rG\cdot \rW)_{t})\\[0.3cm]
\hspace*{6.5cm}
+\dfrac{\partial }{\partial x}\big (\varrho 
\dfrac{\partial (\rG\cdot \rW)_{t}}{\partial x}\big ) \\[0.3cm]
y(\cdot ,0)=u_{0}(\cdot ).
\end{array}%
\right.   
\label{eq. Equivalent determinist}
\end{equation}
\par
In this way, the study is reduced to find sufficient conditions on which
this deterministic problem is well posed. 
Let us present the main result for the simpler case of $(\rE_{\beta, \rG })$ for which we suppose, additionally, 
\begin{equation}
g:\RR\rightarrow \RR\text{ is continuous increasing and }%
g(r)r\geq \mathrm{C}_{g}r^{2}\text{ for some }
\rC_{g}>0.  \label{eq:Hypo g}
\end{equation}
For the Sellers model we assume 
\begin{equation}
\beta :\RR\rightarrow \RR \text{ is a continuous nondecreasing
	and bounded function.}  \label{eq:Hypo beta Sellers}
\end{equation}
Then we have
\begin{theo}[The stochastic Sellers model]
\label{theo:mainSellers} Assume \eqref{eq:infinitestochasticnoise},\eqref{eq:Hypo Wiener cylindrical},
\eqref{eq:infinitestochasticnoiseG},\eqref{eq:Hypo W}, \eqref{eq:paraisometria}, 
\eqref{eq:Hypo beta Sellers}, \eqref{eq:Hypo f(t,x)}, \eqref{eq:Gpredictable},  \eqref{eq:Hypo g},  as well as
\begin{equation}
f\in \rL^{2}(0,\rT;\rV'). 
\label{eq:Hypo f(t,x)V}
\end{equation}
Then, for each $u_{0}\in \rL^{2}(\rI),$ there exists a unique analytically weak solution $u\in \cC([0,\rT];\rH)\cap \rL^{2}(0,\rT;\rV)$ of $(\rE_{\beta, \rG }).$ If, in addition,
\begin{equation}
f\in \rL^{2}(0,\rT;\rH), 
\label{eq:Hypo H f}
\end{equation}
\begin{equation}
\dfrac{\partial }{\partial x}\left( \varrho \dfrac{\partial }{\partial x}
(\rG\cdot\rW)_{t}\right)\in \rL^{2}(0,\rT;\rH),  
\label{eq:Hypo xderivate W}
\end{equation}
then $u\in \cC([0,\rT];\rH)\cap \rL^{2}(0,\rT;\rV)$, $u$ is a mild solution and $y\in \rH^{1}(\delta ,\rT;\rH)$ for any $\delta \in (0,T)$ (with $y$ defined by \eqref{eq:Change of variables}). Moreover, in this last case, if $\widehat{u}\in \cC([0,\rT];\rH)\cap \rL^{2}(0,\rT;\rV)$ is the
corresponding solution to $\widehat{u}_{0}\in \rL^{2}(\rI)$
and $\widehat{f}\in \rL^{2}(0,\rT;\rH)$ then, for any $t\in [0,\rT],$
\begin{equation}
\begin{array}{l}
\disp \left\Vert u_{t}-\widehat{u}_{t} \right\Vert _{\rL
	^{\infty }(0,\rT:\rL^{2}(\rI))}\leq e^{t\rQ\rS_{o}\rL_{\beta }}\left( \left\Vert  u_{0}-\widehat{u}_{0} \right\Vert _{\rL^{2}(\rI)}+\int_{0}^{t}
\left\Vert  f_{s}-\widehat{f}_{s}\right\Vert _{\rL^{2}(\rI)}ds\right) ,\\ [.35cm]
\disp \left\Vert \left[ u_{t}-\widehat{u}_{t}\right] _{+}\right\Vert _{\rL%
	^{\infty }(0,\rT:\rL^{2}(\rI))}\leq e^{t\rQ\rS_{o}\rL_{\beta }}\left( \left\Vert \left[ u_{0}-\widehat{u}_{0}\right] _{+}\right\Vert _{\rL^{2}(\rI)}+\int_{0}^{t}
\left\Vert \left[ f_{s}-\widehat{f}_{s}\right] _{+}\right\Vert _{\rL^{2}(\rI)}ds\right) , 
\end{array} 
\label{eq:comparisonestimate}
\end{equation}%
where $\rL_{\beta }$ is the Lipschitz constant of function $\beta $,
and $\left[ \cdot \right] _{+}$ denotes the positive part of a function.
\end{theo}
The proof will be given in Section \ref{Sec:Sellers} as well as
some remarks on possible improvements. Notice that the last of the above estimates allows to conclude that 
$$
u_{0}(x)\leq \widehat{u}_{0}(x)\text{ and }f_{t}(x)\leq \widehat{f}_{t}(x)%
\text{ implies that }u_{t}\leq \widehat{u}_{t}\text{ for } a.e.~\omega \in
\Omega .
$$
Our quantitative comparison estimate \eqref{eq:comparisonestimate} improves
the comparison result obtained in \cite{BucPar} for the case of Dirichlet
boundary conditions and a "space-time white noise" ($\rW_{x,t}$
standard Brownian sheet to which it can be associated the fundamental
solution $\rG_{t}(x,y)$ of the heat equation). We notice that in the special case of the finite dimensional stochastic noise the existence of mild solutions only requires to assume that $\Phi (\cdot,t)\in\rD(\rA)^{m}$.

The treatment of the Budyko case is more delicate and will be presented in
Section \ref{Sec:Budyko}. Finally, in Section \ref{Sec:applications} we will
give two applications. We will get a different version (more
general in some sense) of the previous results by the authors, \cite{DD}, in
which the stochastic term was supposed of the form $\varepsilon {\rG}\dfrac{\partial \rW}{\partial t}$ and the stability of solutions,
when $\varepsilon \rightarrow 0$, was analyzed. It was shown there that the
associated solution $u^{\varepsilon ,b}$ converges to a \emph{solution} $%
u^{b}$ of the deterministic problem. Finally, we will consider
 the asymptotic behaviour of solutions when $
t\rightarrow +\infty $. We recall that the usual way to study the
"stationary solution" associated to a stochastic partial differential
equation\ is by means of the notion of {\em invariant
	measure} associated to the {\em transition
	group $\rP_{t}$ corresponding to a given measure preserving group of
	transformations} in $\Omega,~\theta _{t}:\Omega
\rightarrow \Omega $, such that the map $\left( t,\omega \right) \mapsto
\theta _{t}\omega $ is measurable and satisfies 
\begin{equation*}
\theta _{t+s}=\theta _{t}\circ \theta _{s}=\theta _{s}\circ \theta
_{t},\qquad \theta _{0}=\mathrm{Id.}
\end{equation*}%
Notice that 
\begin{equation*}
\EE\big [\displaystyle{\Vert \mathrm{W}(\cdot ,t)\Vert }_{\mathrm{H}}^{2}%
\big ]=t\sum_{n\geq 1}\dfrac{1}{\mu _{n}}\rightarrow \infty ,\quad \hbox{as
	$t\rightarrow +\infty$},
\end{equation*}%
and thus the abstract formulation (\ref{eq:SDEabs}) must be understood as a
non-autonomous equation with a perturbation which, in some sense, grows as $%
t\rightarrow +\infty $ (see, $e.g.$, \cite{DaPZer} and the
recent monograph~\cite{Liu} with its many references). Nevertheless, we will
prove a result which does not appear to have been observed previously: if%
\textrm{\ }$\lim_{t\rightarrow \infty }\rG(x,t)=0$ it is possible to give a
meaning to the limit $(\rG\cdot \rW)_{\infty }\doteq\lim_{t\rightarrow +\infty
}(\rG\cdot \rW)_{t}$, even if $\theta _{t}:\Omega \rightarrow \Omega $ is the
identity, and then it is possible to characterize the  $\omega -$limit set
as in the case of similar deterministic energy balance models (see, $e.g.$ 
\cite{Diaz-Her-Tello}). When $\lim_{t\rightarrow \infty }\rG(x,t)\neq 0$ (for
instance when $\rG_{t}$ is constant in time) the convergence of the
stochastic term $(\rG\cdot \rW)_{t}$, as $t\rightarrow
+\infty $, has no sense, and hence the study (for a {\em general group of
	transformations} in $\Omega $, $\theta _{t}:\Omega
\rightarrow \Omega $) must follows some conceptually different approaches:
either by the mentioned theory of invariant measures or by the more general
study of the {\em global random attractor} such as in \cite
{Diaz Langa Valero}, for some energy balance models.

The organization of this paper is the following: section 2 will be devoted
to the proof of Theorem~\ref{theo:mainSellers}. The treatment of the case of 
$\beta $ multivalued (the Budyko problem) will be presented in Section 3.
Section 4 will be devoted to develop the above mentioned two applications. Finally in an Appendix we recall different results dealing with cylindrical Wiener processes and the stochastic convolution.

\section{The stochastic Sellers model}

\label{Sec:Sellers}
In contrast with the exposition made in the Introduction, to treat the stochastic Sellers problem it is easier to make a direct treatment by working with the full operator
\begin{equation}
\cA_{\beta }u\doteq -\dfrac{\partial }{\partial x}\left( \varrho 
\dfrac{\partial }{\partial x}u\right) +g(u)-\rQ\rS\beta (u).
\label{eq:Sellersoper}
\end{equation}%
Then the deterministic problem \eqref{eq:SDEabs} can be reformulated as the  time-dependent Cauchy problem 
\begin{equation}
\left\{ 
\begin{array}{l}
\dfrac{d}{dt}y_{t}+\cA_{\beta }(t)y_{t}=f_{t},\quad t>0, \\[0.15cm] 
y_{t}\big |_{t=0}=u_{0}\in \rH,
\end{array}
\right. 
\label{eq:timedependentCauchyproblem}
\end{equation}
depending on a random parameter, $\omega \in \Omega$, where the time-dependent operator is $%
\cA(t)v\doteq \cA_{\beta }(v+(\rG\cdot\rW)_{t})$.
Evaluated on each $\omega \in \Omega $, the function 
$$
t\mapsto y_{t}(\omega )\in \rV
$$
satisfies the (deterministic) Cauchy problem 
$$
\left\{ 
\begin{array}{l}
\dfrac{d}{dt}y_{t}(\omega )+\cA(t)(\omega )y_{t}(\omega
)=f_{t},\quad t>0,\text{ in }\rV^{\prime }, \\[0.15cm] 
y_{t}(\omega )\big |_{t=0}=u_{0}\in \rH.
\end{array}
\right.
$$
Here $\cA(t)(\omega )y_{t}(\omega )=\cA_{\beta }(y_{t}(\omega
)+(\rG\cdot\rW)_{t}(\omega ))\in \rH$. We recall that 
$\rW_{t}$ has not bounded variation although it is H\"{o}lder
continuous for any of exponent $\alpha \in (0,1/2)$ (see, $e.g.$, \cite{DaPZ}).
However the operator $\cA (t)$ is measurable in a simple sense.

The proof of the first part of  Theorem \ref{theo:mainSellers} will be consequence of the
following abstract result (see \cite[Theorem 4.17]{Bar}) which only requires
the measurability in $t$ of the operator under consideration $\cA (t)$,  as well as a
classical variational structural set of assumptions

\begin{theo}[\protect\cite{Bar}]
\label{theo:Thm Barbu 4.17} Let $\{\cA(t);t\in [0,\rT]\}$
with $\cA(t):\rV\rightarrow \rV^{\prime }$ be a family of
nonlinear operators such that

\begin{enumerate}
\item[i)] the function $t\rightarrow \cA(t)u(t)$ is measurable from $%
[0,\rT]$ to $\rV^{\prime }$ for every measurable function $%
u:[0,T]\rightarrow \rV$,

\item[ii)] there exists a $\delta \geq 0$ such that, for $a.e. ~ t\in (0,
\rT), ~\omega \in \Omega $, the operator $u\rightarrow \delta u+
\cA(t,\omega )u$ is monotone and demicontinuous (that is,
strongly-weakly continuous from $\rV$ to $\rV^{\prime }$),

\item[iii)] there exists $\gamma >0$ and $\rC_{1}\geq 0$ such that 
$$
\pe{\cA(t,\omega )u}{u} \geq \gamma \n{u}_{\rV}^{2}+\rC_{1},\quad \forall u\in \rV,t\in [0,\rT],
$$

\item[iv)] there exists $\rC_{2}>0$ such that 
$$
\n{\cA(t)u}_{\rV^{\prime }}\le \rC_{2}\left( 1+\n{u} _{\rV}\right) \quad \forall u\in \rV
,t\in [0,\rT].
$$
\end{enumerate}
Then, for any $y_{0}\in \rH$ and any $f\in\rL^{2}\left( 0,\rT;\rV^{\prime }\right)$, there exists a unique absolutely continuous function $y\in \rW^{1,2}\left( 0,\rT; \rV'\right)$ such that $y\in \cC([0,\rT];\rH)\cap \rL^{2}(0,\rT;\rV)$ and 
$$
\left\{ 
\begin{array}{lc}
\dfrac{dy}{dt}(t)+\cA(t)y(t)=f(t) & a.e.\text{ }t\in (0,\rT) 
\text{, in }\rV^{\prime }, \\[.2cm] 
y(0)=y_{0}. & 
\end{array}%
\right.
$$
\end{theo}
We point out that the proof of Theorem \ref{theo:Thm Barbu 4.17} was obtained in \cite[Theorem 4.17]{Bar} by means of the addition of two maximal monotone operators from $\cV=\rL^{2}(0,\rT; \rV)$ on $\mathcal{V}^{\prime }=\rL^{2}\left( 0,\rT;\rV^{\prime }\right) $ and that one of them is the operator associated to $%
\dfrac{d}{dt}y_{t}(\omega )$.

\par
\noindent 
{\sc Proof of Theorem \ref{theo:mainSellers}}
As indicated before, to prove the existence and uniqueness of an analytically weak solution it is enough to apply the change of variable $y_{t}=u_{t}-(\rG\cdot\rW)_{t}$ (which is completely well justified thanks to the assumptions on $(\rG\cdot\rW)$) and to prove that if we define 
$$
\cA(t)y=\cA_{\beta }(y+(\rG\cdot\rW)_{t})
$$
then all the conditions of Theorem \ref{theo:Thm Barbu 4.17} are fulfilled
once we take 
\begin{equation}
\delta \geq \rQ\rS_{0}\rL_{\beta }.
\label{eq:Hypo delta}
\end{equation}
We introduce the time dependent operator $\cA(t)u\doteq 
\cA_{\beta }(u+(\rG\cdot\rW)_{t})$ (with $\cA%
_{\beta } $ given by \eqref{eq:Sellersoper}). Certainly 
$$
t\rightarrow \cA(t)u(t)\quad 
\hbox{is measurable for every measurable function 
	$u:[0,\rT]\rightarrow \rV$},
$$
so that i) holds. Moreover, from the monotonicity of $g$ we have 
$$
\begin{array}{ll}
\displaystyle{\langle \cA u-\cA v,u-v\rangle 
}_{\rV^{\prime }\times \rV} & \hspace*{-0.2cm}\displaystyle%
=\int_{\rI}\varrho \big (u_{x}-v_{x}\big )^{2}dx \\[0.2cm] 
& +\displaystyle\int_{\rI}\big (g(u+\big (\rG\cdot\rW
\big )_{t})-g(v+(\rG\cdot\rW)_{t})\big )\big (u-v\big )dx \\[0.3cm] 
& \displaystyle+\int_{\rI}\rQ\rS(x)\big (\beta
(u+(\rG\cdot\rW)_{t})-\beta (v+(\rG\cdot\rW%
)_{t})\big )(u-v)dx \\[0.3cm] 
& \hspace*{-0.2cm}\displaystyle\geq \int_{\rI}\varrho \big (%
u_{x}-v_{x}\big )^{2}dx-\rQ\rS_{0}\rL_{\beta }\big )%
\int_{\rI}|u-v|^{2}dx%
\end{array}
$$
therefore 
$$
\displaystyle{\langle \cA u+\delta u-\cA v-\delta v,u-v\rangle }_{\rV^{\prime }\times \rV}\geq \int_{%
	\rI}\varrho \big (u_{x}-v_{x}\big )^{2}dx+\big (\delta -\rQ 
S_{0}\rL_{\beta }\big )\int_{\rI}|u-v|^{2}dx\geq 0,
$$
provided \eqref{eq:Hypo delta}, which proves the monotonicity of $\cA(t)u+\delta u$ (condition ii)).
\par
\noindent On the other hand, concerning the coercivity condition iii) we point out
that by introducing 
$$
\left\{ 
\begin{array}{l}
\widehat{f}_{t}\doteq \rQ\rS\rM+f_{t}, \\[0.15cm] 
\widehat{\beta }(u)\doteq \beta (u)-\rM, \\[0.15cm] 
\widehat{\cA}_{\beta }u\doteq -\dfrac{\partial }{\partial x}\left(
\varrho \dfrac{\partial }{\partial x}u\right) +g(u)-\rQ\rS
\widehat{\beta }(u),
\end{array}
\right.
$$
the problem under consideration can be equivalently reformulated as 
$$
du+\widehat{\cA}_{\beta }\big (u+(\rG\cdot\rW)_{t}\big )dt=%
\widehat{f}_{t}dt,
$$
so that, it is enough to prove the conditions iii) and iv) for operator $%
\widehat{\cA}_{\beta }(t)u+\delta u.$ Concerning iii), using (H$%
_{\beta }$) we have that $\hat{\beta}((\rG\cdot\rW)_{t})\leq
0$ and since $\rS\geq \rS_{0}$, we get 
$$
\int_{I}\rQ\rS\hat{\beta}(u+(\rG\cdot \rW)_{t})(u+(\rG\cdot\rW)_{t})dx\le \rQ\rS_{0}
\rL_{\beta }\int_{\rI}|u|^{2}dx. 
$$
Then, by assumption \eqref{eq:Hypo g} 
$$
\begin{array}{ll}
\displaystyle{\langle \cA_{\beta }u,u\rangle }_{\rV^{\prime
	}\times \rV} & \hspace*{-0.2cm}\displaystyle=\int_{\rI}\varrho
u_{x}^{2}dx+\int_{\rI}g(u)udx+\int_{\rI}\rQ\rS%
\beta (u)udx \\[0.3cm] 
& \hspace*{-0.2cm}\displaystyle\geq \int_{\rI}\varrho u_{x}^{2}dx+(%
\mathrm{C}_{g}+\delta -\rQ\rS_{0}\rL_{\beta })\int_{%
	\rI}|u|^{2}dx,%
\end{array}
$$
and thus 
$$
\displaystyle{\langle \widehat{\cA}_{\beta }u+\delta u,u\rangle }_{
	\rV'\times \rV}\geq \min \{1,\rC_{g}+\delta -%
\rQ\rS_{0}\rL_{\beta }\}\left( \int_{\rI}\varrho
u_{x}^{2}dx+\int_{\rI}|u|^{2}dx\right) 
$$
provided 
$$
\rC_{g}+\delta \geq \rQ\rS_{0}\rL_{\beta }. 
$$
Finally, to prove iv) let $v$ such that $\Vert v\Vert _{\rV}=1$. Then 
$$
\begin{array}{ll}
|\langle {\widehat{\cA}_{\beta }(t)u},v\rangle _{V^{^{\prime
	}}V}|&\hspace*{-.2cm}\disp \leq \left\vert \int_{\rI}\varrho (x) \frac{\partial u}{\partial x}
\big (u+(\rG\cdot\rW)_{t}\big )\right\vert +\left\vert \int_{
	\rI}g(u+(\rG\cdot\rW)_{t})vdx\right\vert \\ [.35cm]
& \disp+\left\vert
\int_{\rI}\rQ\rS\hat{\beta}(u+(\rG\cdot \rW)_{t}vdx\right\vert . 
\end{array}
$$
Applying H\"{o}lder inequality 
$$
\begin{array}{ll}
\displaystyle\left\vert \int_{\rI}\varrho (x)\frac{\partial }{\partial x}
\big (u+(\rG\cdot\rW)_{t}\big )\frac{\partial v}{ \partial x}
dx\right\vert & \hspace*{-.2cm}\displaystyle\leq \int_{I}\bigg |\sqrt{\varrho }
\frac{\partial u}{\partial x}\bigg |\cdot \bigg |\sqrt{\varrho }\frac{\partial v
}{\partial x}\bigg |dx\\ [.35cm]
&\hspace*{-.2cm}\disp +\int_{I}\bigg |\sqrt{\varrho }\frac{\partial (
	\rG\cdot\rW)_{t}}{\partial x}\bigg |\cdot \bigg |\sqrt{\varrho }\frac{
	\partial v}{\partial x}\bigg |dx \\[0.4cm] 
& \hspace*{-.2cm}\displaystyle\leq \bigg (\left\Vert \frac{\partial u}{%
	\partial x}\right\Vert _{\rL^{2}(\rI:\varrho )}+\left\Vert (\rG\cdot\rW)_{t}\right\Vert _{\rV}\bigg )\Vert v_{x}\Vert _{
	\rL^{2}(\rI:\varrho )} \\[0.4cm] 
& \hspace*{-.2cm}\leq \big (\Vert u\Vert _{\rV}+\rC\big )\Vert
v\Vert _{\rV},\qquad t\in \lbrack 0,\rT].
\end{array}
$$
On the other hand, since the problem is one-dimensional, $\rV\subset 
\rL^{\infty }(\rI)$, and then the mere continuity of $g$ and
the assumption \eqref{eq:Hypo W} imply that 
$$
\left\vert \int_{\rI}g(u+(\rG\cdot\rW%
)_{t})vdx\right\vert \leq \big (\Vert u\Vert _{\rV}+\mathrm{C}\big )%
\Vert v\Vert _{\rV},\qquad t\in [0,\rT]. 
$$
Analogously
$$
\left\vert \int_{\rI}\rQ\rS\hat{\beta}(u+(\rG\cdot
\rW)_{t})vdx\right\vert \leq \rQ\rS_{1}\Vert \hat{%
	\beta}\Vert _{\infty }\int_{I}|v|dx\leq \rQ\rS_{1}\Vert
_{\infty }\Vert \hat{\beta}\Vert _{\infty }\sqrt{2}, 
$$
so that, condition iv) holds. 
\par
\noindent Thanks to the assumptions on $(\rG\cdot\rW)$, when $f\in \rL^{2}(0,\rT;\rH)$, the analytically weak solution is a mild solution (see \cite[Appendix A]{LiRo}) and the regularity $y\in \rH^{1}((0,\rT];\rH)$ is consequence of Corollary 3.3.2 of \cite{Vra} when the problem is
reformulated in terms of 
\begin{equation}
\left\{ 
\begin{array}{l}
\dfrac{d}{dt}y_{t}-\dfrac{\partial }{\partial x}\left( \varrho \dfrac{%
	\partial }{\partial x}y_{t}\right) =\rF_{t}(y_{t})\quad t>0, \\[0.15cm] 
y_{0}=u_{0}\in \rH,%
\end{array}%
\right. \quad \omega \in \Omega,  \label{eq:Pb Sellers via Vrabie}
\end{equation}
with 
$$
\rF_{t}(y_{t})=\rQ\rS\beta (y_{t}+(\rG\cdot\rW)_{t})-g(y_{t}+(\rG\cdot\rW)_{t})+f_{t}. 
$$
In order to prove the quantitative comparison estimate \eqref{eq:comparisonestimate} we observe that%
\begin{equation}
\left\{ 
\begin{array}{l}
\dfrac{d}{dt}(y_{t}-\widehat{y}_{t})+\left[ \rA_{\beta
}(y_{t}+(\rG\cdot\rW)_{t})-\rA_{\beta }(\widehat{y}%
_{t}+(\rG\cdot\rW)_{t})\right] =f_{t}-\widehat{f}_{t}\quad
t>0, \\[0.15cm] 
y_{0}-\widehat{y}_{0}=u_{0}-\widehat{u}_{0}\in \rH,%
\end{array}%
\right.  \label{eq:Eq difference y}
\end{equation}
(recall the change of notation introduced in \eqref{eq:leadingpartHoperator}). Thus 
$$
\begin{array}{ll}
\rA_{\beta }(y_{t}+(\rG\cdot\rW)_{t})-\rA%
_{\beta }(\widehat{y}_{t}+(\rG\cdot\rW)_{t}) & \hspace*{%
	-0.2cm}\displaystyle=-\dfrac{\partial }{\partial x}\left( \varrho (x)\dfrac{%
	\partial }{\partial x}(y_{t}-\widehat{y}_{t})\right)\\ [.4cm]
& \hspace*{-.7cm} +\left[ g(y_{t}+(
\rG\cdot\rW)_{t})-g(\widehat{y}_{t}+(\rG\cdot\rW)_{t})\right] \\[0.3cm] 
& \hspace*{-.7cm} -\rQ\rS(x)\left[ \beta (y_{t}+(\rG\cdot\rW)_{t})-\beta (\widehat{y}_{t}+(\rG\cdot\rW)_{t})
\right] .
\end{array}
$$
Multiplying \eqref{eq:Eq difference y} by $(y_{t}-\widehat{y}_{t})_{+}$ and
arguing as in the proof of the monotonicity we arrive to 
$$
\frac{1}{2}\frac{d}{dt}\left\Vert \left[ u_{t}-\widehat{u}_{t}\right]
_{+}\right\Vert _{\rL^{2}(\rI))}^{2}\leq \rQ\rS%
_{0}\rL_{\beta }\left\Vert \left[ u_{t}-\widehat{u}_{t}\right]
_{+}\right\Vert _{\rL^{2}(\rI)}+\left\Vert \left[ f_{t}-%
\widehat{f}_{t}\right] _{+}\right\Vert _{\rL^{2}(\rI)}.
$$
Then, by Gronwall's inequality we arrive to the conclusion. The proof of the other inequality is similar but multiplying now \eqref{eq:Eq difference y} by $y_{t}-\widehat{y}_{t}.\fin$

\begin{rem}\rm A stronger version of the assumptions \eqref{eq:Hypo W} and \eqref{eq:Hypo xderivate W} are imposed in the paper \cite{BGuessLRock} for other purposes. Notice that assumption \eqref{eq:Hypo xderivate W} must be understood as a restriction on the $x$-dependence of the term $\rG(t,x)$ involved at the noise. For instance,
this condition is trivially satisfied in the case of a finite dimensional
stochastic noise of the form $\rG=\rG(t)$ with $\rG\in \rL^{\infty }(0,\rT)$ and $%
\rW(x,t)=\sum_{i=1}^{m}\rB^{i}(t)\be_{i}(x)$ with $\rB^{i}(t)\in \rL^{2}(0,\rT)$.
Indeed, it suffices to use that $\be_{i}(x)$ is the corresponding orthonormal Legendre
polinomial (see \eqref{eq:baseA}) and thus $\dfrac{\partial }{\partial x}\left( \varrho \dfrac{%
	\partial \be_{i}}{\partial x}\right) \in \rL^{\infty }(\rI).\fin$ 
\end{rem}

\begin{rem}\rm Notice that in the formulation \eqref{eq:Pb Sellers via Vrabie} it
is enough to have a $\rL^{2}$-time dependence of the term $\rF_{t}(r)$ and that, by the
contrary, a reformulation of the problem using the subdifferential theory, since 
$$
\widehat{\rA}_{\beta }(t)u+\delta u=\partial \phi (t,u),
$$
would require that the dependence of the convex function $\phi (t,u)$ with
respect to the time let of bounded variation (see \cite{Yotsutani} and its
references), which does not occur since $\rW_{t}$ is merely H\"{o}lder continuous
in time.$\fin$ 
\end{rem}

\section{On the stochastic Budyko  model}

\label{Sec:Budyko}

Now we consider the Budyko model in which $\beta $ is defined by the multivalued
maximal monotone graph given by \eqref{eq:BetaHeaviside}. Theorem~\ref{theo:Thm Barbu 4.17} can not be directly applied but some different
strategies can be followed in order to prove the existence of a mild solution when working with the $\rH$-realization $\rA:\rD(\rA)\rightarrow \rH$ of the operators, as defined in \eqref{eq:leadingpartHoperator}. One possibility is to adapt to our framework the previous results of \cite{DaPF}.  A more direct line of proof is to approximate $%
\beta $ by a Lipschitz maximal monotone graph $\beta _{\lambda }$. For
instance we can use the Yosida approximation of $\beta ,$ 
\begin{equation}
\beta _{\lambda }=\frac{1}{\lambda }(I-(I+\lambda \beta )^{-1}),
\label{eq:Yosidabeta}
\end{equation}%
where $\rI$ denotes the identity operator. It is easy to see that $\beta
_{\lambda }$ satisfies (H$_{\beta }$) $\ $in the sense that $\beta _{\lambda
}$ is a \emph{bounded} maximal monotone graph in $\RR^{2}$, with $%
\beta _{\lambda }(s)\in \lbrack m,\rM]$, $\forall s\in \RR$
(see, $e.g.$ \cite{Bre}, page 45). 
\begin{rem}\rm It is not too difficult to check that instead to use the Yosida approximation $\beta _{\lambda }$ it is possible to use a more specific approximation of this maximal monotone graph as, for instance, the {\em sigmoidal function} (as it was used, $e.g.$ in \cite{HK1,HK2} for a different problem involving the Heaviside discontinuous function). Nevertheless, the Yosida approximation is applicable for any maximal monotone graph $\beta $  and leads to convergence results in an easier way. $\fin$
\end{rem}
We consider, then, the (unique) solution $
y_{t}^{\lambda }$ of the corresponding approximate problem 
\begin{equation}
\left\{ 
\begin{array}{l}
\dfrac{d}{dt}y_{t}+\rA_{g}(t)y_{t}=\rQ\rS(x)\beta
_{\lambda }(y_{t}+(\rG\cdot\rW)_{t})+f_{t}+\dfrac{\partial 
}{\partial x}\left( \varrho \dfrac{\partial (\rG\cdot \rW)_{t}}{\partial x}\right) ,\quad t>0, \\[0.15cm]
y_{0}=u_{0}\in \rH,%
\end{array}%
\right. 
\label{eq:eqapproximate}
\end{equation}
where now 
$$
\rA_{g}(t)v\doteq-\dfrac{\partial }{\partial x}\left( \varrho \dfrac{%
\partial v}{\partial x}\right) +g\big (v+(\rG\cdot\rW)_{t}\big ).
$$
Arguing as in Lemma 2 of (\cite{Diaz Escorial}) (see also \cite{DD}) we
have
\begin{theo}[The stochastic Budyko model]
Assume \eqref{eq:infinitestochasticnoise}, \eqref{eq:Hypo Wiener cylindrical}
\eqref{eq:infinitestochasticnoiseG}, \eqref{eq:Gpredictable}, \eqref{eq:paraisometria},  \eqref{eq:Hypo g}, \eqref{eq:Hypo H f} and \eqref{eq:Hypo xderivate W}
\label{theo:Budyko approx}
Let $u_{0}\in \rL^{2}(\rI),$ and let $u^{\lambda }\in \cC([0,\rT];%
\rH)\cap \rL^{2}(0,\rT;\rV)$ be the unique mild
solution of \eqref{eq:eqapproximate}. Then there exists $u\in \cC([0,\rT];\rH)\cap \rL^{2}(0,\rT;\rV)$
solution of the stochastic Budyko problem $(\rE_{\beta,\rG }),$ with $\beta $ given
by \eqref{eq:BetaHeaviside}, such that $u^{\lambda }\rightharpoonup u$ in $%
\rL^{2}(0,\rT;\rV)$ and strongly in $\cC([0,\rT];\rH).\fin$
\end{theo}

\bigskip

\begin{rem}\rm 
Under the above conditions, it is possible to apply the
iterative method of super and subsolutions to get the existence of a minimal
and a maximal solutions as in \cite{Diaz Escorial}. As a mater of fact this
method was systematically used in the study of the stochastic problem in 
\cite{DD} for non-cylindrical Wiener processes. $\fin$ 
\end{rem}

\begin{rem}\rm
It is possible to obtain the existence of solutions to problem $(\rE_{\beta,\rG })$ with $\beta $, given by \eqref{eq:BetaHeaviside}, through a
fixed point argument. That was carried out (for \ $\rG\equiv $0$ $ and $%
u_{0}\in \rL^{\infty }$) in \cite{DT} for the more general case in which the
spatial domain is a compact Riemannian manifold without boundary (a study
which includes the case considered in this paper in which the manifold is a
sphere and the solution is only dependent on the latitude). The changes to
consider a non-autonomous problem (corresponding to a stochastic term $\rG\cdot\rW$
satisfying the set of conditions indicated above) are quite standard and will not be presented here.
In fact, problems of this type were considered in Chapter 4 of \cite{Vra}
(see also \cite{DV}), under the general formulation 
$$
\left\{ 
\begin{array}{l}
\dfrac{d}{dt}y_{t}+\rA(t)y_{t}(\omega )\in \rF_{t}(y_{t}),\quad t>0,%
\text{ in }H, \\[0.15cm]
y_{0}(\omega )=u_{0}\in \rH,
\end{array}
\right. 
$$
where $\rA(t)=\partial \varphi (t,u)$ and $\rF_{t}(y_{t})$ is $\varphi 
$-demiclosed. It is not difficult to show (see \cite{Diaz Escorial}) that $\rA_{g}(t)v=\partial \varphi (t,v)$ for a suitable convex function $\varphi $ and that 
$$
\rF_{t}(y_{t})=\rQ\rS(x)\beta \big (y_{t}+(\rG\cdot\rW)_{t}\big )+f_{t}+%
\dfrac{\partial }{\partial x}\left( \varrho (x)\dfrac{\partial (\rG\cdot \rW)_{t}}{\partial x}\right) 
$$
is $\varphi $-demiclosed in the sense of \cite{Vra}, but we will not enter here in
the details.$\fin$ 
\end{rem}
Concerning the uniqueness of mild solutions to the stochastic Budyko problem, since
assumption~\eqref{eq:Hypo W} is perfectly compatible with the condition that 
$\rG$ vanishes locally in a neighborhood of the {\em free boundary}
of the solution (the set $\cF=\left\{ (x,t):u(x,t)=-10\right\} $)
then the concrete examples presented in (\cite{Diaz Escorial} and \cite{DT}) can be
adapted to  show that, in general, there is multiplicity of solutions of the
problem $(\rE_{\beta ,\rG}),$ with $\beta $ given by %
\eqref{eq:BetaHeaviside}. Nevertheless, it is possible to get the uniqueness
of mild solutions in the class of "nondegenerate functions":

\begin{definition}
Let $w\in \rL^{\infty }(\rI)$. We say that $w$ satisfies the strong
(resp. weak) nondegeneracy property if there exists $\rC>0$ and $%
\varepsilon _{0}>0$ such that for any $\varepsilon \in (0,\varepsilon _{0})$ 
$$
\left\vert \{x\in \rI:|w(x)+10|\leq \varepsilon \}\right\vert \le \rC\varepsilon
$$
(\mbox{resp. }$|\{x\in \rI:0<|w(x)+10|\leq \varepsilon \}|\le \rC\varepsilon )$.
\end{definition}

By adapting to the corresponding deterministic problem the results of \cite{Diaz Escorial} and \cite{DT} we
have:

\begin{theo}
\label{theo:Thm uniqueness Budyko} Let $u_{0}\in \rL^{\infty }(%
\rI)$ and assume the same conditions than in Theorem \ref{theo:Budyko approx}.
\begin{enumerate}
	\item[\bf i)] Assume that there exists a mild solution $u(\cdot ,t)$ satisfying the
	strong nondegene\-racy property for any $t\in [0,\rT]$. Then $u$ is
	the unique bounded mild solution of problem $(\rE_{\beta ,\rG%
	}) $ with $\beta $ given by~\eqref{eq:BetaHeaviside}.
	
	\item[\bf ii)] At most there is a unique mild solution of the mentioned problem among the
	class of bounded mild solutions satisfying the weak nondegeneracy property.$\fin$
\end{enumerate}
\end{theo}

The main tool in the proof of Theorem \ref{theo:Thm uniqueness Budyko} \ is the fact that
under the nondegeneracy property the multivalued term generates a
continuous operator from $\rL^{\infty }(\rI)$ into $\rL^{q}(\rI)$, for any $q\in [1,\infty )$.

\begin{lemma}[\protect\cite{Diaz Escorial}]
\label{ldmma:Lemm Budyko} \ 

\begin{enumerate}
\item[\bf i)] Let $w,\hat{w}\in \rL^{\infty }(\rI)$ and assume
that $w$ satisfies the strong nondegeneracy property. Then for any $q\in
[1,\infty )$ there exists $\widetilde{\rC}>0$ such that for any $z,%
\hat{z}\in \rL^{\infty }(\rI),~z(x)\in \beta (w(x)),~\hat{z}
(x)\in \beta (\hat{w}(x))~a.e.~x\in \rI$ we have 
\begin{equation}
\parallel z-\hat{z}\parallel _{\rL^{q}(\rI)}\leq
(a_{f}-a_{i})\min \{\widetilde{\rC}\parallel w-\hat{w}\parallel _{\rL^{\infty }(\rI)}^{1/q}~,~2^{1/q}\}.  
\label{eq:lema3}
\end{equation}

\item[\bf ii)] If $w,\hat{w}\in \rL^{\infty }(\rI)$ satisfy
the weak nondegeneracy property then 
\begin{equation}
\int_{\rI}(z(x)-\hat{z}(x))(w(x)-\hat{w}(x))dx\leq (a_{f}-a_{i})%
\widetilde{\mathrm{C}}\parallel w- \hat{w}\parallel _{\rL^{\infty
	}(\rI)}^{2}.  \label{eq:debil}
\end{equation}
\fineqnum
\end{enumerate}
\end{lemma}
\par
{\sc Proof of Theorem \ref{theo:Thm uniqueness Budyko}} It is enough to
reproduce the proof given in \cite{Diaz Escorial} (see also \cite{DT}) with the same
adaptation to the case $\rG\neq {0}$ than the arguments used
in the proof of the monotonicity of the operator given in the proof of Theorem \ref
{theo:mainSellers}. $\fin$

\begin{rem}\rm Curiously enough, it is possible to prove the uniqueness of a
different (and larger) class of solutions when the stochastic term has a
different nature (see \cite{DD}). $\fin$ 
\end{rem}
\section{Two applications}

\label{Sec:applications}

The non-autonomous equation in the equivalent formulation 
\eqref{eq. Equivalent determinist} to the stochastic problem $(\rE_{\beta ,\rG}) $ allows to extend to the stochastic case (after some generalizations)
some of the results available on deterministic energy balance models. In this section we
collect two of them.

\subsection{Stability to the deterministic models}
\label{subsec:stability}
One of the main goals of the paper \cite{DD} was to give a rigorous proof of
the passing from stochastic to deterministic formulations when the Wiener
process is of the form $\varepsilon \dot{\rW}$ for a {\em
	space-time white noise} $\dot{\rW}$ on a filtered probability space $%
(\Omega ,\cF,\{\cF_{t}\}_{t\geq 0},{\PP}).$ In the
special case in which the stochastic term is given by $\varepsilon \rG
d\rW_{t}$ with $\rG\in \cL_{2}(\rH)$
acting on a cylindrical Wiener process $\rW$ satisfying %
\eqref{eq:Hypo W} the mentioned convergence of solutions, as $\varepsilon
\rightarrow 0$, is a simple consequence of the abstract results on the
convergence of maximal monotone operators. Let us explain this for the case
of the Sellers stochastic model. For a given $\varepsilon >0$, let 
\begin{equation}
\rA_{\beta }^{\varepsilon }(t)v\doteq -\dfrac{\partial }{\partial x}%
\left( \varrho \dfrac{\partial }{\partial x}v\right) +g(v+\varepsilon \big (\rG\cdot \rW)_{t}\big )-\rQ\rS\beta (v+\varepsilon
(\rG\cdot\rW)_{t})-\varepsilon \dfrac{\partial }{\partial x}%
\left( \varrho \dfrac{\partial (\rG\cdot\rW)_{t}}{\partial x%
}\right) ,
\end{equation}%
for $v\in \rD(\rA_{\beta }^{\varepsilon }(t))$ and $a.e.~t\in (0,%
\rT)$, where $\rA_{\beta }^{\varepsilon }$ is the associated maximal monotone operator, from $\rH$ to $\rH
$, with 
$$
\rD(\rA_{\beta }^{\varepsilon }(t))=\big \{ v\in \rV\text{ such
	that }\rA_{\beta }^{\varepsilon }(t)v\in \rH\big\} \quad a.e.~t\in (0,\rT).
$$
We will apply the abstract result given in \cite{BeCrPa} (Theorem 11.1,
Proposition 11.2 and E11.17). See also the pioneering Th\'{e}or\`{e}me 3.16 of \cite{Bre} on the convergence of autonomous operators.

\begin{theo}[\protect\cite{BeCrPa}]
\label{theo:Thm Conver BCP} Let $\rA_{n}(t),~t\in [0,\rT%
],~n=1,2,...,\infty $ be a family of operators such that $\delta I+\rA%
_{n}(t)$ are m-accretive in a Banach space $\mathrm{X}$ and denote by $%
J_{n,\lambda }^{\delta ,t}$ the resolvent operators of $\delta I+\rA%
_{n}(t),~i.e.$ 
$$
J_{n,\lambda }^{\delta ,t}=(I+\lambda (\delta I+\rA_{n}(t)))^{-1}.
$$
Assume that there exists $h\in \cC([0,\rT])$ and a
nondecreasing function $L:[0,+\infty )\rightarrow [0,+\infty )$ such that 
\begin{equation}
\left\Vert J_{n,\lambda }^{\delta ,t}v-J_{n,\lambda }^{\delta
	,s}v\right\Vert _{\rX}\leq \lambda \left\Vert h(t)-h(s)\right\Vert
L(\left\Vert v\right\Vert )\text{ for }t,s\in [0,\rT],v\in \mathrm{X}.
\label{eq:Hypo grow J}
\end{equation}
Assume also that 
\begin{equation}
\lim_{n\rightarrow \infty }J_{n,\lambda }^{\delta ,t}v=J_{\infty ,\lambda
}^{\delta ,t}v\text{ for any }v\in \rD,  \label{eq:Hypo converg resolve}
\end{equation}%
where $\rD$ is a dense set of $\mathrm{X}.$ Let $u_{n}$ be the mild solution
of the problem 
$$
\left\{ 
\begin{array}{l}
\dfrac{du_{n}}{dt}+\rA_{n}(t)u_{n}\ni f_{n}(t)\quad t>0,\text{ in }%
\mathrm{X}, \\[0.15cm] 
u_{n}(0)=u_{0,n},
\end{array}
\right.
$$
and let $u_{\infty}$ be the mild solution of 
$$
\left\{ 
\begin{array}{l}
\dfrac{du_{\infty }}{dt}+\rA_{\infty }(t)u_{\infty }\ni f_{\infty}(t)\quad t>0,\text{ in }\rX, \\[0.15cm] 
u_{\infty }(0)=u_{0,\infty },%
\end{array}%
\right.
$$
with $f_{n},f_{\infty }\in \rL^{1}(0,\rT:\mathrm{X})$ and $%
u_{0,n}\in \overline{\rD(\rA_{n}(0))},~u_{0,\infty }\in \overline{\rD(%
	\rA_{\infty }(0))}.$ Asumme, finally, that 
\begin{equation*}
f_{n}\rightarrow f_{\infty }\text{ in }\rL^{1}(0,\rT:\rX)\text{
	and }u_{0,n}\rightarrow u_{0,\infty }\text{ in }\rX.
\end{equation*}
Then $u_{n}\rightarrow u_{\infty }$\textit{\ in }$\cC([0,\rT]:\rX).\fin$
\end{theo}

\bigskip

As a consequence we have

\begin{theo}
\label{:theoThm convergencia sellers} Assume the conditions of Theorem \ref{theo:mainSellers}. 
For $n\in N$ let $u^{n}$ be the solution of problem $(\rE_{\beta ,\frac{1}{n}\rG})$. Then, when $n \rightarrow \infty$, it follows $u^{n}\rightarrow u^{\infty}$ in $\cC([0,\rT]:\rH)$, where $u^{\infty }$ is the solution of the deterministic problem $(\rE_{\beta ,0})$.
\end{theo}

{\sc Proof} Let $\varepsilon =\frac{1}{n}$. Due to the
peculiar expression of operators $\rA_{\beta }^{\varepsilon }(t)$ we
know that $\rD(\rA_{\beta }^{\varepsilon }(t))=\rD(\rA_{\beta
}^{0}(t))$ is independent on $t$. Moreover, since $(\rG\cdot\rW)_{t}$ is H%
\"{o}lder continuous we also know that condition \eqref{eq:Hypo grow
J} holds. Thus, to get the conclusion it is enough to apply Theorem \ref%
{theo:Thm Conver BCP} once we check that the resolvent convergence %
\eqref{eq:Hypo converg resolve} is fulfilled. Let $v\in \rH$ and let $%
w^{\varepsilon }$ be the unique solution of the stationary problem in $\rI$
\begin{equation}
-\dfrac{\partial }{\partial x}\left( \varrho \dfrac{\partial }{\partial x}%
w^{\varepsilon }\right) +g(w^{\varepsilon }+\varepsilon (\rG\cdot
\rW)_{t})-\rQ\rS\beta (w^{\varepsilon }+\varepsilon (\rG\cdot \rW)_{t}+\delta w^{\varepsilon }=v+\varepsilon \dfrac{%
\partial }{\partial x}\left( \varrho \dfrac{\partial (\rG\cdot \rW)_{t}}{\partial x}\right) 
\label{eq:Eq Stationary resolvent eps}
\end{equation}
and let $w^{\infty }$ be the unique solution of the equation 
$$
-\dfrac{\partial }{\partial x}\left( \varrho \dfrac{\partial }{\partial x}%
w^{\infty }\right) +g(w^{\infty })-\rQ\rS\beta (w^{\infty
})+\delta w^{\infty }=v,\text{ in }\mathrm{I.}
$$
Then by multiplying equation \eqref{eq:Eq Stationary resolvent eps} by $
w^{\varepsilon }$ and using the monotonicity and coercivity of $\delta I+%
\rA_{\beta }^{\varepsilon }(t)$ we get the uniform estimate 
$$
\left\Vert w^{\varepsilon }\right\Vert _{\rV}\leq \rC ,
$$
for some $\rC>0$. Thus $w^{\varepsilon }\rightharpoonup w$ in $
\rV$ to some $w\in \rV$ and then also strongly in $\rH$
and in $\rL^{\infty }(\rI)$. In addition, $\beta
(w^{\varepsilon }+\varepsilon (\rG\cdot\rW)_{t})$ is
bounded (from assumption (H$_{\beta }$)) and since $g$ is continuous and $%
(\rG\cdot\rW)_{t}$ is bounded we get that 
$$
g(w^{\varepsilon }+\varepsilon (\rG\cdot\rW)_{t})\text{ is
bounded in }\rI.
$$
Since $g$ and $\beta $ are maximal monotone operators they are closed for
the weak-weak convergence and thus $g(w^{\varepsilon }+\varepsilon (\rG\cdot\rW)_{t})\rightharpoonup g(w)$ and $\beta
(w^{\varepsilon }+\varepsilon (\rG\cdot\rW)_{t})\rightharpoonup \beta (w)$ a.e. in $\rI$, and then strongly in $\rH$. So that, by the uniqueness of solutions $w=w^{\infty }$ and
condition (\ref{eq:Hypo converg resolve}) holds.$\fin$
\par
To prove a similar conclusion for the case of the Budyko type
model we will follow a different strategy:

\begin{theo}
\label{theo:Thm convergencia Budyko} Assume the conditions of Theorem \ref
{theo:Budyko approx}. For $n\in N$ let $u^{n}$ be any solution of problem $(
\rE_{\beta ,\frac{1}{n}\rG})$. Let $f\in \rL^{2}(0,%
\rT;\rH)$ and assume (\ref{eq:Hypo xderivate W}). Then $
u^{n}\rightarrow u^{\infty }$ in $\cC([0,\rT]:\rH)$,
as $n\rightarrow \infty $, where $u^{\infty }$ is a solution of the
deterministic problem $(\rE_{\beta ,0})$ with $\beta $ given
by \eqref{eq:BetaHeaviside}.
\end{theo}
\par
\noindent 
{\sc Proof} Let $\varepsilon =\frac{1}{n}$. Let $%
u^{\varepsilon }$ be any solution of problem $(\rE_{\beta
,\varepsilon \rG})$ we use the notation $u_{t}^{\varepsilon
}=w_{t}^{\varepsilon }+\varepsilon(\rG\cdot\rW)_{t}$ and 
$$
\xi _{t}^{\varepsilon }=-g(w^{\varepsilon }+\varepsilon (\rG\cdot \rW)_{t})+\rQ\rS b^{\varepsilon }+\varepsilon \dfrac{%
\partial }{\partial x}\left( \varrho \dfrac{\partial (\rG\cdot \rW)_{t}}{\partial x}\right) +f_{t},
$$
where $b^{\varepsilon }(x,t)\in \beta (w^{\varepsilon }(x,t)+\varepsilon
(\rG\cdot\rW)_{t})$, for a.e. $x\in \rI$ and $t\in
(0,\rT)$, is the section making the identity in equation $(\rE%
_{\beta ,\varepsilon \rG})$. Then $\xi ^{\varepsilon }$ is uniformly
bounded 
\begin{equation*}
\left\Vert \xi ^{\varepsilon }\right\Vert _{\rL^{2}(0,\rT:%
\rL^{2}(\rI))}\leq \mathrm{C},
\end{equation*}%
for some $\mathrm{C}>0$ independently of $\varepsilon $. Thus, there exists $
\xi \in \rL^{2}(0,\rT:\rL^{2}(\rI))$ and a
subsequence (denoted again by $\xi ^{\varepsilon }$) such that $\xi
^{\varepsilon }\rightharpoonup \xi $ in $\rL^{2}(0,\rT:\rL^{2}(\rI))$. As in the Proposition 1 of \cite{Diaz Escorial} the pure diffusion part (see \eqref{eq:leadingpartHoperator})
$$
\rA v=-\dfrac{\partial }{\partial x}\left( \varrho \dfrac{\partial v%
}{\partial x}\right) ,\qquad \rD(\rA)=\left\{ v\in \rV:%
\rA v\in \rL^{2}(\rI)\right\}, \qquad 
$$
can be written as the subdifferential $\rA v =\partial \varphi (v)$
of the convex and lower semicontinuous functional, $\varphi :\rL^{2}(%
\rI)\rightarrow \RR\cup \{+\infty \},$ given by 
\begin{equation}
\varphi (w)=\left\{ 
\begin{array}{ll}
\displaystyle\frac{1}{2}\int_{\rI}\varrho (x)\left\vert \dfrac{\partial v%
}{\partial x}\right\vert ^{2}dx & \mbox{if }v\in \rV, \\ [.35cm]
+\infty  & \mbox{otherwise}.%
\end{array}%
\right.   \label{eq:defphi}
\end{equation}
Moreover, by Lemma 1 of \cite{Diaz Escorial}, $\partial \varphi (w)$ generates a 
\emph{compact} semigroup of contractions on $\rL^{2}(\rI)$.
Thus, since, for $\varepsilon =\frac{1}{n},$ 
$$
\left\{ 
\begin{array}{l}
\dfrac{du_{n}}{dt}+\partial \varphi (u^{n})\ni \xi ^{n}(t)\quad t>0,\text{
in }\rH, \\[0.15cm]
u^{n}(0)=u_{0},
\end{array}%
\right. 
$$
by Corollary 2.3.2 of \cite{Vra} we know that $u^{n}\rightarrow u^{\infty }$
in $\cC([0,\rT]:\rH)$ with $u^{\infty }$ solution of 
$$
\left\{ 
\begin{array}{l}
\dfrac{du^{\infty }}{dt}+\partial \varphi (u^{\infty })\ni \xi _{t}\quad t>0,%
\text{ in }\rH, \\[0.2cm]
u_{0}^{\infty }=u_{0}.%
\end{array}%
\right. 
$$
But $u^{n}\rightarrow u^{\infty }$ in $\cC([0,\rT]:\rH)
$, as $n\rightarrow \infty $, and the facts that $\beta $ is a maximal
monotone graph and $g$ is continuous imply that $\xi _{t}(x)\in
-g(u_{t}^{\infty }(x))+\rQ\rS(x)\beta (u_{t}^{\infty
}(x))+f_{t}(x),$ for $a.e.~x\in \rI$ and $t\in (0,\rT).$ Thus $%
u^{\infty }$ is a solution of the deterministic problem $(\rE_{\beta ,0})$ with $\beta $ given by \eqref{eq:BetaHeaviside} and the proof ends.$\fin$

\begin{rem}\rm We recall that in \cite{Diaz Escorial} it was obtained some regularity results on the solution of the auxiliary linear problem with pure diffusion: for any $u_{0}\in\rL^{2}(\rI)$ there exists a unique function $u\in\cC([0,\rT]:\rL^{2}(\rI))$, such that $u(t)\in\rD(\rA)$ for $a.e. t\in (0,\rT] $, $t^{\frac{1}{2}}\dfrac{d u(t)}{dt}\in \rL^{2}\big (0,\rT:\rL^{2}(\rI)\big )$ and is a mild solution of
\begin{equation}
\left\{ 
\begin{array}{l}
\dfrac{d u}{d t}(t)+\rA u(t)=0,\quad t>0,\\[0.2cm] 
u(0)=u_{0}\in\rH.
\end{array}
\right. 
\label{eq:deterSellerproblemleading}
\end{equation}	
Moreover, if $u_{0}\in\rL^{p}(\rI),~1\le p\le \infty$ then $u(t)\in\rL^{p}(\rI)$ for $a.e.~ t\in (0,\rT] $. In fact, if $u_{0}\in\rV$ one has $\dfrac{d u}{dt}\in \rL^{2}\big (0,\rT:\rL^{2}(\rI)\big )$. Finally, the application 
$\rS(t)u_{0}=u(t)$ defines a compact semigroup of contractions on $\rL^{2}(\rI).\fin$
\label{rem:semigroup}
\end{rem}
\subsection{On the asymptotic behaviour of solutions as $t\rightarrow
	+\infty $}
A pedagogical way for studying the behaviour of solutions $u$ of the stochastic
problem $(\rE_{\beta ,\rG}),$ as $t\rightarrow +\infty $, consists in to inquire into this question previously in the framework of deterministic problems. So we first consider the deterministic problem 

\begin{equation}
\hspace*{-.2cm}
\left\{ 
\begin{array}{ll}
\dfrac{d}{dt}y_{t}-\dfrac{\partial }{\partial x}\left( \varrho \dfrac{
\partial y_{t}}{\partial x}\right) +g(y_{t}+(\rG\cdot \rW
)_{t})&\hspace*{-.2cm}\disp \in \rQ\rS\beta (y_{t}+(\rG\cdot \rW)_{t})+f_{t}\\
&\quad\disp +
\dfrac{\partial }{ \partial x}\left( \varrho \dfrac{\partial (\rG\cdot\rW)_{t}}{\partial x}\right) , \\[0.15cm] 
y_{0}=u_{0}\in \rH, &
\label{eq:SDErandom}
\end{array}%
\right.
\end{equation}%
under similar conditions to the ones of the paper \cite{Diaz-Her-Tello}. It concerns a special class
non-autonomous dynamical systems whose asymptotic limits (in time) are autonomous differential equations. For the moment, we are not assuming that $(\rG\cdot \rW)$ is related to any  cylindrical Wiener process but are simply some given datum. We need some extra requirements on the data. Here we will assume for a while

\begin{description}
\item[(H$_{f,\rW}$)] $f\in \rL^{\infty }(\rI\times (0,\infty ))$ and 
\begin{equation}
\dfrac{\partial }{\partial x}\left( \varrho \dfrac{\partial (\rG\cdot \rW)_{t}}{ \partial x}\right)\in \rL^{\infty }(\rI\times (0,\infty )),
\end{equation}
\item[(H$_{\infty }$)] there exists $f_{\infty }\in \rV^{\prime }$
and $(\rG\cdot \rW)_{\infty}\in \rV$ with 
\begin{equation}
(\rG\cdot\rW)_{\infty}=\lim_{t\rightarrow +\infty}\int_{0}^{t}\rG_{s}d\rW_{s},
\label{eq:integralinfty}
\end{equation}
such that
\begin{equation*}
\dfrac{\partial }{\partial x}\left( \varrho \dfrac{\partial (\rG\cdot \rW)_{\infty}}{\partial x}\right) \in \rL^{\infty }(\rI),
\end{equation*}
$$
\left \{
\begin{array}{ll}
\displaystyle \int_{t-1}^{t+1}\Vert f(\tau ,\cdot )-f_{\infty }(\cdot )\Vert
_{V^{\prime }}d\tau \rightarrow 0 & \quad \mbox{ as }t\rightarrow \infty , \\%
[.275cm] 
\displaystyle \int_{t-1}^{t+1}\left[\left  \Vert \dfrac{\partial }{\partial x}%
\left( \varrho \dfrac{\partial (\rG\cdot\rW)_{t}}{\partial x}%
\right) -\dfrac{\partial }{ \partial x}\left( \varrho \dfrac{\partial (\rG\cdot\rW)_{\infty}}{\partial x}\right) \right \Vert _{\rV^{\prime }}%
\right] d\tau \rightarrow 0 & \quad \mbox{ as }t\rightarrow \infty ,\text{ }
\\[.275cm ] 
\displaystyle \int_{t-1}^{t+1}\Vert (\rG\cdot\rW)_{\tau}-(\rG\cdot\rW)_{\infty}\Vert _{\rV}d\tau \rightarrow 0 & \quad %
\mbox{ as }t\rightarrow \infty .%
\end{array}%
\right .
$$
\end{description}

\begin{rem}\rm  
Condition \eqref{eq:integralinfty} requires some words on its meaning in the framework of cylindrical Wiener processes. Notice than if we assume \eqref{eq:integralinfty} then, by the isometry of Wiener, we get (see the Appendix below)
$$ 
\lim_{t\rightarrow \infty}\EE \left  [\left \Vert (\rG\cdot\rW)_{t}\right \Vert ^{2}_{\HH}\right  ]=\lim_{t\rightarrow \infty}\cQ_{t}^{\rG\cJ}=\cQ_{\infty}^{\rG\cJ},
$$
so that, necessarily, we must have  
\begin{equation}
\cQ_{\infty}^{\rG\cJ}=\int^{\infty}_{0}\n{\rG_{s}\cJ}_{\cL_{2}(\HH)}^{2}ds<+\infty.
\label{eq:covarianceoperatorinfty}
\end{equation}
For instance, in the special case in which  $\rG_{t}=\psi(t)\rG$, for some constant operator $\rG$, one has 
$$
\cQ_{t}^{\rG\cJ}=\n{\rG\cJ}_{\cL_{2}(\HH)}^{2}\int^{\infty}_{0}\big (\psi (s)\big )^{2}ds .
$$
This introduces a constraint which can be satisfied in many cases. For instance, for the choice $\psi(t)=(t+a^{-1})^{-\alpha}$, with $2\alpha >1$, one has
$$
\cQ_{\infty}^{\rG\cJ}=\dfrac{a^{2\alpha-1}}{2\alpha-1}\n{\rG\cJ}_{\cL_{2}(\HH)}^{2},
$$
so that,
$$
\cQ_{\infty}^{\rG\cJ}<+\infty\quad \Leftrightarrow \quad \n{\rG\cJ}_{\cL_{2}(\HH)}^{2}<+\infty.
$$
In this case, the relative Wiener isometry becomes
\begin{equation}
\EE\big[\n{(\rG\cdot\rW)_{\infty}}^{2}_{\rL^{2}(\rI)}\big ]=\dfrac{a^{2\alpha-1}}{2\alpha-1}\n{\rG\cJ}_{\cL_{2}(\HH)}^{2}.
\label{eq:Wienerisometry} 
\end{equation}
This explains that, in this example, $(\rG\cdot \rW)_{\infty }\neq 0$ even if $\psi (t)\rightarrow 0
$ as $t\rightarrow +\infty $ (the important convergence is the one of $(\rG\cdot \rW)_{t}$ and not the pointwise convergence of $\rG_{t}$). In other words: $(\mathrm{G}\cdot \mathrm{W})_{\infty }\neq (\rG
_{\infty }\cdot \rW)_{\infty }$ if $\rG_{\infty
}(x)=\lim_{t\rightarrow \infty }\rG(x,t)=0$. Notice that, in fact, if $\rG_{\infty }(x)=\lim_{t\rightarrow \infty }\rG(x,t)\neq 0$ then $(\rG
_{\infty }\cdot \rW)_{\infty }$ has no sense.  We will see later a different study of the stability in terms of random attractors for the case in which $\rG_{t}$ is constant in time.  $\fin$
\label{rem:asympt hipo noise}
\end{rem}

As in Lemma 1 of \cite{Diaz-Her-Tello}, we can obtain the global regularity of the solutions on $(0,\infty )$ when the data, and in particular the time derivative of $(\rG\cdot\rW)_{t}$, is in $\rL^{1}(%
\rI)$: something that in general fails in the framework of stochastic equations but for which the property \eqref{eq:solucionhasatinfinito} below can be obtained by other ways. More precisely

\begin{prop}
\label{Proposi global existence} Assume $u_{0}\ \in \rV\cap \rL^{\infty
}(\rI),~  f\in \rL^{\infty }(\rI\times (0,\infty ))\cap
W_{loc}^{1,1}(0,\infty ;\rL^{1}(\rI))$ \mbox{ and } 
\begin{equation}
\hspace*{-.3cm}
\left \{ 
\begin{array}{l}
\displaystyle\int_{t}^{t+1}\left \Vert {\frac{\partial f}{\partial t}}(s,\cdot
)\right \Vert _{\rL^{1}(\rI)}ds\leq \rC_{f},\ \forall t>0\\
\mbox{ where }\rC_{f}
\mbox{
is a time independent constant, } \\[.45cm ] 
(\rG\cdot\rW) \in \rL^{\infty }(0,\infty ;\rV)\cap \rW_{loc}^{1,1}((0,\infty );\rL^{1}(\rI)), \\[.3cm] 
\displaystyle\int_{t}^{t+1}\left[ \bigg \Vert \dfrac{\partial ^{2}}{\partial
t\partial x}\left( \varrho \dfrac{\partial (\rG\cdot\rW)_{t}%
}{\partial x}\right) \bigg \Vert _{\rL^{1}(\rI)}\right] ds \leq
\rC_{(\rG\cdot \rW)},\ \forall t>0\\
\mbox{ where the constant }\rC_{(\rG\cdot \rW)}\mbox{
is independent on $t$. }%
\end{array}
\right .
\label{eq:hypo time derivability}
\end{equation}
Then there exists a solution of \eqref{eq:SDErandom} verifying 
\begin{equation}
y\in \rL^{\infty }(0,\infty ;\rV)\ \mbox{ and }\ \frac{\partial y}{
\partial t}\in \rL^{2}(0,\infty ;\rL^{2}(\rI)).
\label{eq:solucionhasatinfinito}
\end{equation}
\fineqnum
\end{prop}
The following theorem (which is a direct application of Theorem 1 of \cite{Diaz-Her-Tello}) proves
the stabilization of solutions $y$ satisfying the deterministic problem \eqref{eq:SDErandom} in the class of bounded functions given by \eqref{eq:solucionhasatinfinito}.  We
define the $\omega $-limit set of $y$ by 
$$
\omega (y)=\{y_{\infty }\in \rV\cap \rL^{\infty }(\rI):\ \exists
t_{n}\rightarrow +\infty \mbox{ such that }y(t_{n},\cdot )\rightarrow
y_{\infty }\mbox{ in }\rL^{2}(\rI)\}.
$$
\begin{theo}\label{theo:infinito}
Let $u_{0}\in \rL^{\infty }(\rI)\cap \rV$, let $u$ be a mild solution of of the stochastic problem $(\rE_{\beta ,\rG}),$  and let $y$ be the corresponding bounded solution of \eqref{eq:SDErandom} 
satisfying \eqref{eq:solucionhasatinfinito}. Then
\begin{itemize}
	\item[\bf i)] $\omega (y)\neq \emptyset $,
	\item[\bf ii)] \textit{if }$y_{\infty }\in \omega (u)$ then $\exists 
	t_{n}\rightarrow +\infty $ such that  $y(t_{n}+s,\cdot )\rightarrow
	y_{\infty }$ in $\rL^{2}(-1,1;\rL^{2}(\rI))$ and $y_{\infty
	}\in $\textit{\ }$\rV$. Moreover if we define $u_{\infty }=y_{\infty
	}+(\rG\cdot\rW)_{\infty}$ then $u_{\infty }$ is a bounded solution of the
	deterministic stationary equation of the asymptotic autonomous problem,  
	$$
	(\rP_{\rQ})\qquad -\dfrac{\partial }{\partial x}\left( \varrho \dfrac{\partial
		u_{\infty }}{\partial x}\right) +g(u_{\infty })\in \rQ\rS\beta (u_{\infty
	})+f_{\infty }\mbox{  in }\rI.
	$$
	\item[\bf iii)] \textit{In fact, if }$y_{\infty }\in \omega (y)$\textit{\ then }
	$\exists \{\hat{t}_{n}\}\rightarrow +\infty $\textit{\ such that }$y(\hat{t}
	_{n},\cdot )\rightarrow y_{\infty }$ strongly in$\ \rV.\fin$
\end{itemize}
\end{theo}
As mentioned in the Introduction, the usual treatment of the stochastic model (avoiding to assume that the
time derivative of $(\rG\cdot\rW)_{t}$ is in $\rL^{1}(\rI))$ and assumption \eqref{eq:integralinfty}) 
requires the application either of the theory of invariant measures or of the random dynamical systems theory. This approach was already developed in \cite{Diaz Langa Valero}. Here we will extend the results of such paper to the more general formulation of problem $(\rE_{\beta , \rG})$ such as it was presented in the Introduction.
In the rest of this section we weal deal merely with a special type of finite dimensional stochastic noise. 
Although in our case we will work with the Hilbert space $\rH$ the abstract theory deals with a general complete and separable metric space ($\rX,d_{\rX})$  with the Borel $\sigma $-algebra $\cB\left( \rX\right) $. We consider now a probability space $(\Omega ,\cF,%
\PP)$ and we suppose given a measure preserving group of transformations in $\Omega $, $\theta _{t}:\Omega \rightarrow \Omega $,  such that the map $\left( t,\omega \right) \mapsto \theta _{t}\omega $ is measurable and satisfies
$$
\theta _{t+s}=\theta _{t}\circ \theta _{s}=\theta _{s}\circ \theta
_{t} \hbox{ and }\theta _{0}=\hbox{\rm Id}.
$$
In order to define later the notion of attractor, in what follows, time variable $t$ takes values in $\RR$ endowed with the Borel $\sigma $-algebra $\cB( \RR).$ We denote by $\cP(\rX) $ the set of non-empty (non-empty closed) subsets of $\rX.$ 
The type of finite dimensional stochastic noise we will consider now is 
\begin{equation}
\rG(x)\dfrac{\partial \rW(x,t)}{\partial t}\doteq\Phi(x)\dfrac{d \rB(t)}{dt},
\quad (x,t)\in\rI\times (0,\infty),
\label{eq:finitestochasticnoise2}
\end{equation}
where $\Phi \in\rD(\rA)$ (see \eqref{eq:leadingpartHoperator})
and with  $\rB(t)$ a real standard Brownian motion on the filtration space $\{\Omega ,\cF, \{\cF_{t}\}_{t\ge 0}, \PP \}$. As a matter of fact, we must assume now that the Wiener process is "two-sided", because $t$ takes values in $\RR$.
\par
The main idea to define the random dynamical system is to consider the map
$\cG:\RR_{+}\times \Omega \times 
\rX\rightarrow \rX $ such that to an initial datum $u_{0}\in \rH$, $\omega\in\Omega$ and $t>0$ makes correspond one solution $u(t,\cdot)$ of the stochastic problem $(\rE_{\beta , \rG})$ particularized in the parameter $\theta _{t}(\omega)$. Nevertheless, since in the case of the stochastic Budyko problem we may have multiplicity of solutions (recall Section 3), the definition of random dynamical system must be enlarged to treat with multivalued functions. More in general we introduce the following notion.
\begin{definition}
\label{def:MRDS} A set valued map $\cG:\RR_{+}\times \Omega \times 
\rX\rightarrow \cP(\rX)$ is called a Multivalued Random Dynamical System (\rm MRDS in short) if it is measurable, i.e., if given $x\in \rX$ the map  
\begin{equation*}
\left( t,\omega ,y\right) \in \Omega \mapsto dist(x,\cG(t,\omega ,y))\text{ }
\end{equation*}
is measurable, where $\disp dist(x,\rA)=\inf_{d\in \rA}d_{\rX}(x,d)$, for $\rA\subset \rX,$ 
and satisfies
\begin{enumerate}
	\item[\bf i)] $\cG(0,\omega,\cdot)=\hbox{\rm Id}$ on $\rX$;
	
	\item[\bf ii)] $\cG(t+s,\omega ,x)=\cG(t,\theta _{s}\omega ,\cG(s,\omega,x)),\quad \forall t,s\in \RR^{+},~x\in \rX,~\PP-a.s.$ \quad (cocycle property).
	
\end{enumerate}
\end{definition}

We will present later a result indicating that the above conditions hold in the above mentioned special case associated to solutions of $(\rE_{\beta , \rG})$, but before will recall the notion of {\em attractor} set for a general
stochastic partial differential (see the monograph \cite{Carvallo-Langa-Robinson}) and specially the case in which there is a a multivalued term (reason why sometime the multivalued equation is called as {\em inclusion}). 
Some previous definitions are needed:
\begin{definition}
$\cG(t,\omega ,\cdot)$ is said to be upper
semicontinuous if for all $t\in \RR_{+}$ and $\PP-a.s,~  
\omega \in \Omega $, given $x\in \rX$ and a neighbourhood of $\cG(t,\omega ,x)
, ~\cO(\cG(t,\omega ,x)),$ there exists $\delta >0$ such that if $d_{\rX}(x,y)<\delta $ then  
$$
\cG(t,\omega ,y)\subset \cO(\cG(t,\omega ,x)).
$$
On the other hand, $\cG(t,\omega ,\cdot)$ is called lower semicontinuous if for all $t\in \RR_{+}$ and $\cP-a.s~\omega \in \Omega ,$ given $x_{n}\rightarrow x$ $(n\rightarrow +\infty )$ 
and $y\in \cG(t,\omega ,x),$ there exists $y_{n}\in \cG(t,\omega ,x_{n})$ such 
that $y_{n}\rightarrow y.$ It is said to be \textit{continuous} if it is upper and lower semicontinuous.
\end{definition}

Let us denote an $\varepsilon$-neighborhood of a set $\rA$ by $\cO
_{\varepsilon}(\rA) =\left\{ y\in \rX:dist\left( y,\rA\right)
<\varepsilon\right\} $.

\begin{definition}
$\cG(t,\omega,\cdot)$ is said to be \textit{upper
	semicontinuous if in the definition of }upper semicontinuity we replace the 
neighborhood $\cO$ by an $\varepsilon$-neighborhood $\cO_{\varepsilon}$.
\end{definition}

It is clear that any upper semicontinous map is $\varepsilon$-upper
semicontinuous. The converse is true if $G$ has compact values.

\par
Following Crauel and Flandoli \cite{Crauel1} we  introduce now the generalization of the concept of random attractor to the case of a Multivalued Random Dynamical System.  We will recall later a general result for the existence and uniqueness of attractors.  Firstly we need some other definitions.

\begin{definition}
A closed random set $\rD$ is a measurable map $\rD:\Omega \rightarrow  
\cP(\rX) $.
A closed random set $\rD(\omega)$ is said to be  negatively
(resp. strictly) invariant for the {\rm MRDS} if  
$$
\rD(\theta_{t}\omega)\subset \cG(t,\omega)\rD(\omega)\ (\text{resp. }\rD(\theta
_{t}\omega)=\cG(t,\omega)\rD(\omega))\quad \forall\ t\in\RR_{+},~\PP-a.s.
$$
\end{definition}
For convenience we are using the notation $\cG(t,\omega)x\doteq\cG(t,\omega,x)$
Suppose the following conditions on the Multivalued Random Dynamical System $\cG$:

\begin{itemize}
\item[({\bf H1})] There exists an absorbing random compact set $\rB(\omega)$, that 
is, for every bounded set $\rD\subset \rX$, there exists $t_{\rD}(\omega)$ such 
that for all $t\geq t_{\rD}(\omega)$ 
\begin{equation}
\cG(t,\theta_{-t}\omega)\rD\subset \rB(\omega)\quad \PP-a.s.
\end{equation}

\item[({\bf H2})] $\cG(t,\omega):\rX\rightarrow \cC(\rX)$ is upper semicontinuous, for all
$t\in\RR_{+}$ and $\PP-a.a.$ $\omega\in\Omega$.
\end{itemize}

Define the limit set of a bounded set $\rD\subset \rX$ as 
\begin{equation}
\disp \Lambda_{\rD}(\omega)=\bigcap_{\rT\geq0}\overline{\bigcup_{t\geq
\rT}\cG(t,\theta_{-t}\omega)\rD}.
\label{eq:Langa}
\end{equation}
The following results are well-known (see, $e.g.$, \cite{K}):

\begin{lemma}
$\Lambda _{\rD}(\omega )$ is the set of limits of all converging sequences $
\left\{ x_{n}\right\} _{n\geq 1}$, where $x_{n}\in \cG(t_{n},\theta 
_{-t_{n}}\omega )\rD$ as $t_{n}\nearrow +\infty .\fin$
\end{lemma}

\begin{prop}
\label{OmegaLimit}Assume conditions {\bf (H1)} and {\bf (H2)} hold. Then,

\begin{itemize}
\item[\bf i)] $\Lambda _{\rD}( \omega)$ is a non empty compact  subset of $\rB(x)$.

\item[\bf ii)] $\Lambda _{D}(\omega )$ is negatively invariant, that is, $
\cG(t,\omega )\Lambda _{D}(\omega )\supseteq \Lambda _{D}(\theta _{t}\omega )$
for all $t\in \RR^{+},\,\PP-a.s.$ If $\cG(t,\omega )$ is lower 
semicontinuous, then $\Lambda _{D}(\omega )$ is strictly invariant.

\item[\bf iii)] $\Lambda _{\rD}(\omega )$ attracts $\rD,$ that is, $\PP-a.s.$ 
$$
\lim_{t\rightarrow +\infty }\text{dist}(\cG(t,\theta _{-t}\omega )\rD,\Lambda
_{\rD}(\omega ))=0.
$$
\fineq
\end{itemize}
\end{prop}
Finally, we arrive to the mentioned general definition:
\begin{definition}
\label{rand-att} A closed random set $\omega \mapsto  
\cA(\omega)$ is said to be a Global Random Attractor of $\cG$ if:

\begin{itemize}
	\item[\bf i)] $\cG(t,\omega )  \cA (\omega) \supseteq  \cA (\theta_t\omega)
	,$ for all $t\geq 0,~\PP-a.s$ (that is, it is negatively  invariant);
	
	\item[\bf ii)] for all $\rD\subset \rX$ bounded,  
	\begin{equation*}
	\lim_{t\rightarrow +\infty }\text{dist}(\cG(t,\theta _{-t}\omega )\rD, 
	\cA (\omega) )=0;
	\end{equation*}
	
	\item[\bf iii)] $\cA( \omega) $ is compact $\PP-a.s.$ 
\end{itemize}
\end{definition}

We point out that in the previous literature the above type of convergence is usually called as "pullback convergence" and that the corresponding attractors are often called "pullback random attractors".  The following result characterize the random attractor of the Multivalued Random Dynamical Systems (MRDS) (see~\cite{CaLaVa}).

\begin{theo}
\label{theo:ExistAtr}Let assumptions {\bf  (H1)-(H2)} hold, the map $(t,\omega
)\mapsto \overline{\cG(t,\omega)\rD}$ be measurable for all deterministic
bounded sets $\rD\subset \rX$, and the map $x\in \rX\mapsto \cG(t,\omega)x$ have
compact values.  Then, 
\begin{equation}
\cA(\omega)\doteq\overline{\bigcup_{\underset{_{\text{bounded}}}{\rD\subset \rX}}\Lambda
	_{\rD}(\omega)}  \label{eq:globalrandomattractor}
\end{equation}
is a global random attractor for $\cG$ (measurable with respect to $\cF $). It is unique and the minimal closed attracting set. Moreover, if the map $x\longmapsto \cG(t,\omega )x$ is lower semicontinuous  for each fixed $\left( t,\omega \right) $, then the global random attractor $
\cA\left( \omega \right) $ is strictly invariant, i.e., $\cG(t,\omega
) \cA(\omega)= \cA (\theta_t\omega)$ for all $t\geq 0.\fin$
\end{theo}
We point out that there is a conceptual difference between Theorems \ref{theo:infinito} and \ref{theo:ExistAtr}. Although we take a limit when t does to infinity in both cases, in the second case the attraction is when the initial time tends to minus infinity (which is the best we know how to do in the stochastic case). In other words, the dynamics described in the omega limit set before Theorem \ref{theo:infinito} and in \eqref{eq:Langa} are very different. In any case, the key result, which avoid the assumption of regularity that the
time derivative of $(\rG\cdot\rW)_{t}$ is in $\rL^{1}(\rI)$,  will come through the application of the following abstract result (see Theorem \ref{theo:ExistAtr} in  Kapustyan \cite{K}):

\begin{theo}
\label{teho:ExistAtr2}Let $\Omega$ be a metrizable topological space, $\cF
$ be the Borel $\sigma$-algebra, and let $\cG$ satisfying conditions i), ii) of 
Definition \ref{def:MRDS}. Assume the following conditions:

\begin{itemize}
\item[\rm ({\bf H1b})] There exists a measurable mapping $r:\Omega 
\rightarrow\RR_{+}$ such that for $\PP-$almost all $\omega 
\in\Omega$ and for $\rR>0$, there is a $\rT=\rT\left(\rR,\omega\right) >1$ such that
$$\left\Vert \cG( t-1,\theta_{-t}\omega)\bB_{\rR} \right\Vert _{+}\leq 
r\left( \omega\right) $$ 
for all $t\geq \rT$, where $\disp \left\Vert \rA\right\Vert 
_{+}=\sup_{y\in \rA}\left\Vert y\right\Vert $, for $\rA\subset \rX;$
\item[\rm ({\bf H2b})] If $x_{n}\rightarrow x_{0}$ weakly, $
t_{n}\rightarrow t_{0}>0$, $\omega_{n}\rightarrow\omega_{0}$, and $y_{n}\in 
G( t_{n},\omega_{n})x_{n}$, then up to a subsequence $
y_{n}\rightarrow y_{0}\in \cG( t_{0},\omega_{0})x_{0}.$ 
\end{itemize}

\noindent Then $\cG$ generates a Multivalued Random Dynamical Systems and the set $\cA(\omega)
\doteq \overline{\bigcup_{\rR>0}\Lambda_{\bB_{\rR}}(\omega)}$ is a global random 
attractor. It is unique and the minimal closed random attracting set.$\fin$
\end{theo}

We will see that Theorem \ref{teho:ExistAtr2} can be applied to our problem  $(\rE_{\beta ,\rG})$. We note that ii) was carried out in the paper \cite{Diaz Langa Valero} on the special case in which the diffusion operator is replaced by the $1d$-Laplacian operator (and adding Neumann type boundary conditions), $g(r)=\mathrm{C} r$ with $\rC$ a positive constant,  $f(x,t)=f(x)$ and the cylindrical Wiener process as in   \eqref{eq:finitestochasticnoise2}.

Concerning the spatial diffusion, the following result (partially mentioned before) proves that we can replace it by the case with a degenerate weight:

\begin{prop}[\cite{Diaz Escorial}]
\label{prop:Level} The subdiferential $\rA\doteq\partial \varphi $ of the convex function $\varphi$ defined by \eqref{eq:defphi} generates a compact semigroup.$\fin$
\end{prop}

On the other hand, the change of variable \eqref{eq:Change of variables} leads to the deterministic formulation depending on a random parameter presented in previous sections. It is easy to check that defining the multivalued function 
$f:(-1,1)\times\R\rightarrow2^{\RR}$ by
$$
f(x,r)=\rQ\rS(x)\beta(r)-g(r),\text{a.e. }x\in(-1,1),\,r\in\RR.
$$
and the multivalued operator $\rF:\rH\rightarrow\cP({\rH})$,
$$
\rF(u)=\{y+h:y\in \rH, y(x)\in F(x,u(x))~a.e. ~x\in\left(
-1,1\right)  \}
$$
then we have:

\begin{enumerate}
\item[(F1)] $\rF:\rH\rightarrow \cC_{v}(\rH),$ where $\cC_{v}(\rH)$ is the set of nonempty, bounded, closed and convex subset of $\rH$ 
\item[(F2)] There exist $\rC_{1},\rC_{2}\geq0$ such that  $\Vert y\Vert \leq \rC_{1}+\rC_{2}\Vert v\Vert ,$
for all $y\in \rF(v)$, $v\in \rH.$
\item[(F3)] $\rF$ is upper $\epsilon$-semicontinuous in $\rH$.
\item[(F4)] There exist  $\delta>0$ y $\rM>0$ such that
$$
(y,u)\leq-\delta\Vert u\Vert ^{2}+\rM,
$$
for all  $u\in \rD(\rA)$ y $y\in \rF(u).$
\item[(F5)]The level sets $\rM_{\rR}=\left\{  u\in \rD(\varphi)  :\left\Vert u\right\Vert \leq \rR,\varphi(u)  \leq \rR\right\}  $ are compact in $\rH$ for all  $\rR>0.$
\end{enumerate}

As a matter of fact, in \cite{K} it was replaced (F3) by the stronger one that $\rF$ is upper semicontinuous, but in  \cite{Diaz Langa Valero}  it was proved that the generalization given in (F3) is enough to our purposes and that it holds for the version of the problem  $(\rE_{\beta ,\rG})$ considered in that paper. 
\par

Let us show that ${\rF}$ is fine even for a general function $g(u)$. First, it is clear that the growth is sublinear. Indeed, from  (F2) we have 
\begin{equation}
\left\Vert y\right\Vert \leq \rC_{1}+\rC_{2}\left\Vert u\right\Vert +\rC_{2}
\varepsilon\left\Vert \phi\right\Vert \left\vert w\left(  t\right)
\right\vert +\varepsilon\left\Vert A\phi\right\Vert \left\vert w\left(
t\right)  \right\vert \leq\widetilde{\rC}_{1}+\rC_{2}\left\Vert u\right\Vert ,
\label{IneqFNew}
\end{equation}
for all  $y\in\bF(u)  $ and $t\in [0,\rT]$,
where $\widetilde{\rC}_{1}$ depend on  $\omega$ and  $\rT$. Moreover, one deduces from (F1) and (F3) that  $\rF(u)  \in \cC_{v}(\rH)$ and that $u\rightarrow\rF(u)$ is upper $\epsilon$-semicontinuous. Finally, one proves  that, fixed $\omega\in\Omega, $for all $u\in \rH$ the function  $t\mapsto\rF(u(t))$ has a measurable selection.

By using well known results (see, $e.g.$, \cite{Vra}), it is not difficult to show the existence at least a strong solution for each initial data in $\rH$. Furthermore, one may prove that each strong solution can be globally defined for all  
$t\geq0.$

Let $\cD\left(  v_{0},\omega\right)  $ be the set of all strong solution of \eqref{eq:SDErandom}. They are  defined for all $t>0$. Then we construct the multivalued map
$\cG:\RR^{+}\times\Omega\times \rH\rightarrow \cP(\rH)  $ by means of 
$$
\cG\left(  t,\omega,v_{0}\right) \doteq\left\{  v\left(  t\right)  +\varepsilon
\zeta\left(  t\right)  :v\left(  \omega,\text{\textperiodcentered}\right)
\in\mathcal{D}\left(  v_{0},\omega\right)  \right\}  .
$$
One may prove in a similar way as \cite[Proposition 4]{CaLaVa} that  $\cG$ is satisfies the cocycle property (see Definition \ref{def:MRDS}).

Since, by definition,  $\theta_{s}:\Omega\rightarrow\Omega$ is given by 
$\theta_{s}\omega=w\left(  s+\cdot\right)  -w\left(  s\right)  \in\Omega$,
then the function $\widetilde{\zeta}$ relative to $\theta
_{s}\omega$ is
$$
\widetilde{\zeta}\left(  \tau\right)\doteq
\zeta\left(  s+\tau\right)  -\zeta\left(  s\right)  =\phi\left(  w\left(
s+\tau\right)  -w\left(  s\right)  \right).
$$

So that, we arrive to the main result in this section:

\begin{theo}\label{theo:attractor}
\label{AtrClimate} Under the mentioned assumptions, the multivalued random dynamical system $\cG$ associated to problem 
$(\rE_{\beta,\rG})$ has a random global attractor $\cA_{\varepsilon}.\fin$
\end{theo}

Finally, by replacing $\rG$ by $\varepsilon \rG$, as in Section \ref{subsec:stability} it was proved in \cite{Diaz Langa Valero} the following attractor convergence result:

\begin{theo}
We have 
$$
\lim_{\varepsilon\searrow 0}dist(\cA_{\varepsilon}(\omega),
\cA)=0,\quad\forall\omega\in\Omega,
$$
where $\cA$ is the deterministic global attractor corresponding  to $\varepsilon =0.\fin$
\end{theo}

The attractor for $\varepsilon =0$ \ can be better characterized under additional conditions. In particular, motivated by the first part of this subsection, it is relevant to consider the associate stationary problem (independently of conditions \eqref{eq:hypo time derivability} since they are merely sufficient conditions to conclude \eqref{eq:solucionhasatinfinito}). The associate problem, denoted by $(\rP_{\rQ,(\rG\cdot\rW)_{\infty}})$, is given as
$$
-\dfrac{\partial }{\partial x}\left( \varrho (x)\dfrac{%
\partial (y_{\infty }+(\rG\cdot\rW)_{\infty})}{\partial x}%
\right) +g(y_{\infty }+(\rG\cdot\rW)_{\infty})\in \rQ\rS\beta
(y_{\infty }+(\rG\cdot\rW)_{\infty})+f_{\infty }%
\mbox{  in
}\rI.
$$
Notice the this equation coincides with the one of Theorem \ref{theo:infinito}, $(\rP_{\rQ})$ (remember the change of variable). We arrived before to this equation in a different way (see an explicit example in Remark~\ref{rem:asympt hipo noise}).

\par
\noindent  We assume additionally
$$
\hspace*{-4cm}
(\bH_{f_{\infty },(\rG\cdot\rW)_{\infty}}) \qquad f_{\infty },\dfrac{%
\partial }{\partial x}\left( \varrho (x)\dfrac{\partial (\rG\cdot \rW)_{\infty }}{\partial x}\right) \in \rL^{\infty }(\rI)
$$
and that there exist $\rC_{f},\rC_{(\rG\cdot\rW)}>0$ such that 
$$
\left \{
\begin{array}{ll}
-\Vert f_{\infty }\Vert _{\rL^{\infty }(\rI)}\leq f(x)\leq -\rC_{f} & \hspace*{-.4cm}
\text{a.e. }\ x\in \rI \\ [.3cm]
-\left \Vert \dfrac{\partial }{\partial x}\left( \varrho (x)\dfrac{\partial 
(\rG\cdot\rW)_{\infty}}{\partial x}\right)\right \Vert
_{\rL^{\infty }(\rI)}\leq \dfrac{\partial }{\partial x}\left( \varrho 
\dfrac{\partial (\rG\cdot\rW)_{\infty}}{\partial x}\right)
(x)\leq -\rC_{(\rG\cdot\rW)} & \hspace*{-.25cm}\text{a.e. } x\in \rI.
\end{array}
\right .
$$

We also assume some crucial balance among the data:
\begin{itemize}
\item[({\bf H$_{\protect\beta }^{*}$})] $\ \ $there exists two real numbers $0<m<\rM$
and $\epsilon >0$ such that $\beta (r)=\{m\}$ for any$\ \ r\in (-\infty
,-10-\epsilon )$ and $\beta (r)=\{\rM\}$ for any $r\in (-10+\epsilon ,+\infty
) $.

\item[({\bf H$_{\rC_{f}}$})] $g(-10-\epsilon )+\rC_{f}>0$ and $\displaystyle{\frac{%
g(-10+\epsilon )+\Vert f_{\infty }\Vert _{L^{\infty }(\rI)}}{%
g(-10-\epsilon )+\rC_{f}}}\leq {\frac{\rS_{0}\rM}{\rS_{1}m}}.$
\end{itemize}
We recall the notion of solution we will consider: 
\begin{definition} We say that a function $u\in \rV\cap \rL^{\infty }(\rI)$ is a "bounded weak solution" of $(\rP_{\rQ})$ if there exists $z\in \rL^{\infty }(\rI),~z(x)\in \beta (y_{\infty }+(\rG\cdot \rW)_{\infty })$ a.a. $x\in \rI$ such that 
\begin{equation*}
\int_{\rI}\varrho (x)\dfrac{\partial (y_{\infty }+(\rG\cdot \rW)_{\infty })}{\partial x}\dfrac{\partial v}{\partial x}%
dx+\int_{\rI}g(y_{\infty }+(\rG\cdot \rW)_{\infty })vdx=\int_{\rI}\rQ\rS(x)zvdx+\int_{\rI}f_{\infty }vdx,
\end{equation*}%
for any $v\in \mathrm{V}$.
\end{definition}

\par
The following result is, again, a direct application of the
similar result (Theorem~2 proved in~\cite{Diaz-Her-Tello}) and explains the different multiplicity of solutions depending on the parameter $\rQ$:

\begin{theo} \label{th.multi} Let $\bH_{f_{\infty },(\rG\cdot\rW)_{\infty}}$, $(\bH_{\protect\beta }^{*}$) be satisfied. Let $y_{m},y_{M}$be the (unique) bounded solutions of the problems
$$
(\rP_{m})\ \ \ -\dfrac{\partial }{\partial x}\left( \varrho (x)\dfrac{\partial
	(y_{\infty }+(\rG\cdot\rW)_{\infty})}{\partial x}\right)
+g(y_{\infty }+(\rG\cdot\rW)_{\infty})=\rQ\rS m+f_{\infty }\ %
\mbox{ on }\rI,
$$
and
$$
(\rP_{\rM})\ \ \ -\dfrac{\partial }{\partial x}\left( \varrho (x)\dfrac{\partial
	(y_{\infty }+(\rG\cdot\rW)_{\infty})}{\partial x}\right)
+g(y_{\infty }+(\rG\cdot\rW)_{\infty})=\rQ\rS\rM+f_{\infty }\ %
\mbox{ on }\rI,
$$
respectively. Then
\begin{itemize}
\item[i)] for any $\rQ>0$  there is a minimal solution $
\underline{y}$ (resp. a maximal solution $\overline{y}$)
of problem $(\rP_{Q})$. Moreover any other solution $y$ must satisfy 
\begin{equation}
y_{m}\leq\underline{y}\le y\le\overline{y}\ \leq \ y_{\rM}
\label{esti:umM}
\end{equation}
and
$$
\begin{array}{l}
g^{-1}(\rQ\rS_{0}m-\Vert f_{\infty }\Vert _{\rL^{\infty }(\rI)})\leq
u_{m}\leq g^{-1}(\rQ\rS_{1}m-\rC_{f}),  \label{eq:um} \\ [.2cm]
g^{-1}(\rQ\rS_{0}\rM-\Vert f_{\infty }\Vert _{\rL^{\infty }(\rI)})\ \leq
u_{M}\leq g^{-1}(\rQ\rS_{1}\rM-\rC_{f}),  \label{eq:uM}
\end{array}
$$
\item[ii)] for any $\rQ$ there is, at least, a solution $y$ of $(\rP_{\rQ})$ which is a global minimum of the functional
$$
\cJ(w)=\frac{1}{2}\int_{\rI}\varrho (x)\left\vert \dfrac{\partial
	(y_{\infty }+(\rG\cdot\rW)_{\infty})}{\partial x}%
\right\vert ^{2}dx+\int_{\rI}\cG(w+)dA-\int_{\rI}f_{\infty
}wdA-\int_{\rI}\rQ\rS(x)j(w)dA,
$$
on the set  $\rK=\{w\in V,~\cG(w)\in \rL^{1}(I)\}$, where $\beta=\partial j$.
\end{itemize}
\par
\noindent Moreover, if $(\bH_{\rC_{f}})$ holds, and we define the balanced constants
$$
\left \{
\begin{array}{ll}
\disp \rQ_{1}={\frac{\mathcal{G}(-10-\epsilon )+C_{f}}{S_{1}\rM}} &\disp \rQ_{2}={\frac{%
	\mathcal{G}(-10+\epsilon )+\Vert f_{\infty }\Vert _{\rL^{\infty }(\rI)}%
}{S_{0}\rM}}  \label{def.Q} \\ [.4cm]
\disp \rQ_{3}={\frac{\mathcal{G}(-10-\epsilon )+C_{f}}{S_{1}m}} &\disp \rQ_{4}={\frac{%
	\mathcal{G}(-10+\epsilon )+\Vert f_{\infty }\Vert _{\rL^{\infty }(\rI)}%
}{S_{0}m}}
\end{array}
\right. 
$$
we have the following multiplicity of solutions:  
\begin{itemize}
\item[iii)] if $0<\rQ<\rQ_{1}$ then $(\rP_{\rQ})$ has a
unique solution $u_{m}<-10$, which in fact is the
minimum of $\cJ$ on  $\rK$ and in addition
$$
\cG^{-1}(-\Vert f_{\infty }\Vert _{\rL^{\infty }(\rI)})\leq
\lim_{\rQ\searrow 0}\inf \Vert u_{m}\Vert _{\rL^{\infty }(\rI)}\leq
\lim_{\rQ\searrow 0}\sup\Vert u_m\Vert _{L^{\infty }(\rI)}u\leq 
\cG^{-1}(-\rC_{f}),
$$

\item[iv)] if $\rQ_{2}<\rQ<\rQ_{3}$ then $(\rP_{\rQ})$ has
at least three solutions, $u_{i},~i=1,2,3$ with $u_{1}=u_{\rM}>-10$, $u_{2}=u_{m}<-10$ and $u_{1}\geq u_{3}\geq u_{2}$ on $\rI$.
Moreover $u_{1}$ and $u_{2}$ are local minima of $\cJ$ on $\rK\cap \rL^{\infty }(\rI)$,

\item[v)] if $\rQ_{4}<\rQ$ then $(\rP_{\rQ})$ has a
unique solution $u_{\rM}>-10$. Moreover, it is the
minimum of $\cJ$ on $\rK$ and $\Vert u_{\rM}\Vert _{\rL^{\infty }(\rI)}\rightarrow +\infty $ when $\rQ\rightarrow +\infty .\fin$
\end{itemize}
\end{theo}
\begin{rem}\rm  
Sharper multiplicity results require additional details on the data. See, $e.g.$ \cite{ADT, BensidDiaz2, HetzerNato,Hetzer infinite,Hetzer-Paul3,NK,SchmidtBertina}.$\fin$
\end{rem}
\begin{rem}\rm  
At the time of writing this article, the authors do not know if Theorem \ref{theo:ExistAtr} can be generalized to the case of an infinite dimensional noise as considered in Section 2. We also summarize that, for instance, in the example given in Remark \ref{rem:asympt hipo noise} the associated asymptotic problem is purely deterministic and that under suitable conditions there is a different multiplicity of solutions according the values of the parameter $\rQ.\fin$
\end{rem}
\vspace*{-4cm}
\begin{figure}[htp]
\begin{center}
	\includegraphics[width=13cm]{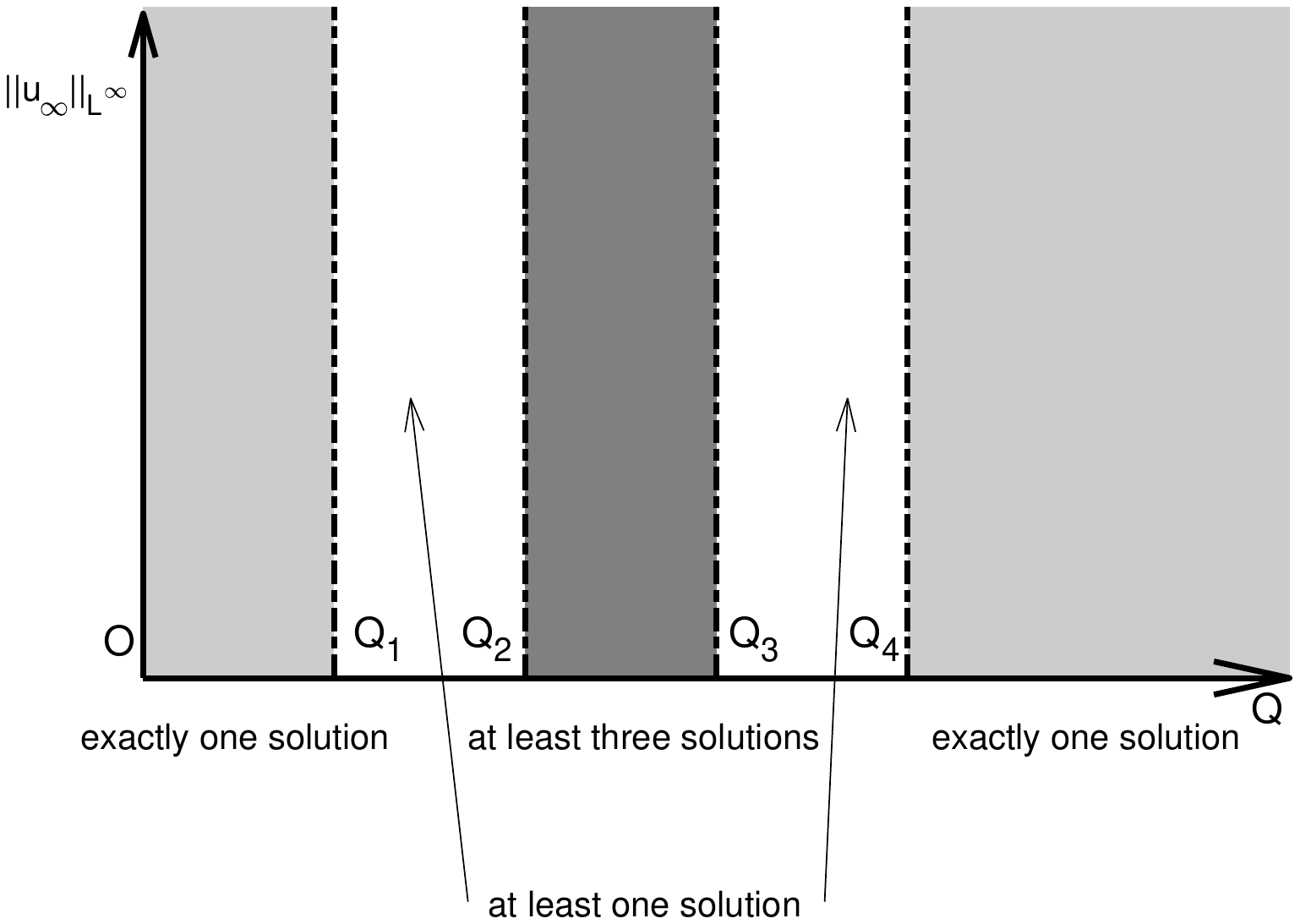}\\ [-6cm]
	\caption{Values of the parameter $\rQ$ with different multiplicity}
\end{center}
\end{figure}

\section{Appendix: On the stochastic framework in this paper}

\subsection{Some few remarks on cylindrical Wiener processes}

\label{Sec:cilindricalWienerprocess}

We consider some cylindrical Wiener process on a separable Hilbert space $\HH_{0}$ by following the definitions and notations of some classical references (see, for instance, \cite{DaPZ}, \cite{Ch} or \cite{LiRo}). So,  $\{\rW_{t}\}_{t\ge 0}$ is given by  
\begin{equation}
\rW_{t}=\sum_{n\ge 1}\rB_{t}^{n}\cJ \widehat{\be}_{n},\quad t\ge 0. 
\label{eq:cylindricalWiener}
\end{equation}
More precisely,  $\cJ:\HH_{0}\rightarrow \HH_{1}$ is a Hilbert-Schmidt embedding, $\cJ\in\cL_{2}(\HH_{0};\HH_{1})$,  involving two separable Hilbert spaces $\HH_{0}$ and $\HH_{1}$ (see Remark  \ref{rem:HIlbertSchamidtTracespace} below), $\{\widehat{\be}_{n}\}_{n\ge 1}$ is a Hilbertian base of $\HH_{0}$ and $\{\rB_{t}^{n}\}_{t\ge 0}$ is a sequence of independent real valued standard Brownian motions.  We note that the property
$$
\begin{array}{ll}
\disp \EE \big  [\big \Vert \sum_{n\ge 1}\rB_{t}^{n}\cJ\widehat{\be}_{n} \big \Vert ^{2}_{\HH_{1}}\big   ]&\disp =\lim_{n\rightarrow \infty}\EE \big  [\big \Vert \sum_{m=1}^{n}\rB_{t}^{m}\cJ\widehat{\be}_{m} \big  \Vert ^{2}_{\HH_{1}}\big  ]\\
&\disp =t\lim_{n\rightarrow \infty}\sum_{m=1}^{n} \n{\cJ\widehat{\be}_{m}}^{2}_{\HH_{1}}=t\n{\cJ}_{\cL_{2}(\HH_{0};\HH_{1})}^{2}<\infty,
\end{array}
$$
shows that for any $t>0$ the series converges in $\rL^{2}(\Omega,\cF,\PP;\HH_{1})$ and uniformly in  $[0,\rT]~\PP-a.s.$ to a Gaussian random variable $\rW_{t}$ with mean $0$ and {\em covariation operator}
$$
\cQ^{\cJ}_{t}=t\n{\cJ}_{\cL_{2}(\HH_{0};\HH_{1})}^{2}.
$$
\begin{rem}\rm \label{rem:HIlbertSchamidtTracespace}Since by definition $\cJ$ is linear and 
$$
\n{\cJ}_{\cL_{2}(\HH_{0};\HH_{1})}^{2}=\sum_{n\ge 1}\left \Vert \cJ \widehat{\be}_{n} \right \Vert ^{2}_{\HH_{1}}=\sum_{n\ge 1}\pe{\cJ \cJ^{*} \widehat{\be}_{n}}{\widehat{\be}_{n}}_{\HH_{1}}=\hbox{trace }\cJ\cJ^{*}<\infty,
$$	
the operator $\cQ^{\cJ}=\cJ\cJ^{*}$ is nonnegative definite and symmetric , with finite trace on $\HH_{1}$. Thus $\cQ^{\cJ}\in \cL_{1}^{+}(\HH_{1})$ and satisfies
$$
\n{\cJ}_{\cL_{2}(\HH_{0};\HH_{1})}^{2}=\hbox{trace }\cJ\cJ^{*}=\n{\cQ^{\cJ}}_{\cL_{1}^{+}(\HH_{1})}.
$$
In fact, the so called {\em Wiener isometry} holds and then
$$
\EE \big  [\big \Vert \sum_{n\ge 1}\rB_{t}^{n}\cJ\widehat{\be}_{n} \big \Vert ^{2}_{\HH_{1}}\big ]=t~\hbox{trace }\cQ^{\cJ}<\infty.
$$
We note that property
\begin{equation}
\EE\big [\pe{\rW_{t}}{u_{1}}_{\HH_{1}}\pe{\rW_{s}}{u_{2}}_{\HH_{1}}\big ]=\min\{s,t\}\pe{\cJ^{*}u_{1}}{\cJ^{*}u_{2}}_{\HH_{0}},\quad u_{1},u_{2}\in\HH_{1},
\label{eq:Winerisometrygeneral}
\end{equation}
also holds. In this framework $\HH_{0}$ is known as the so-called {\em Cameron-Martin} space associated to the centered Gaussian measure related to $\cQ^{\cJ}$.  As it is well known in this theory, the property $\hbox{trace }\cQ^{\cJ}<\infty$ plays a crucial role.$\fin$
\end{rem}
\begin{rem}\rm \label{rem:constructionJ}
In practice, it is very advantageous to start from $\HH_{0}$  and to find a Hilbert space $\HH_{1}$ such that $\HH_{0}$ is densely embedded in $\HH_{1}$ and whose inclusion map
$$
\cJ:~\HH_{0}\rightarrow \HH_{1}
$$
is a Hilbert-Schmidt operator, namely $\cJ\cJ^{*}\in\cL_{1}^{+}(\HH_{1})$. One can always find a Hilbert space $\HH_{1}$ with the above properties. Indeed, let  $\HH_{0}$ be a separable Hilbert space and  $\{\widehat{\be}_{n}\}_{n\ge 1}$ be a Hilbertian base of $\HH_{0}$. For each real sequence $\{\alpha_{n}\}_{n\ge 1}$, with $\disp \sum_{n\ge 1}\alpha_{n}^{2}<+\infty$, one defines the norm 
\begin{equation}
\n{u}_{\HH_{1}}^{2}\doteq\sum_{n\ge 1}\alpha_{n}^{2}|\pe{u}{\widehat{\be}_{n}}_{\HH_{0}}|^{2},\quad u\in\HH_{0},
\label{eq:normH1}
\end{equation}
that clearly satisfies
$$
\n{u}_{\HH_{1}}^{2}\le \n{u}^{2}_{\HH_{0}}\sum_{n\ge 1}\alpha_{n}^{2},\quad u\in\HH_{0}.
$$
In the applications, this norm can be viewed as the one of some negative exponent Sobolev space. 
Then one builds the Hilbert space $\HH_{1}\doteq\overline{\HH_{0}}^{\HH_{1}}$ and  the map $\cJ_{\hbox{\footnotesize Id}}:~(\HH_{0},\n{\cdot}_{\HH_{0}})\rightarrow (\HH_{1},\n{\cdot}_{\HH_{1}})$ given by
setting the identity map $\cJ_{\hbox{\footnotesize Id}} u=u,~u\in \HH_{0}$.
Clearly,  $\{\widehat{\be}_{n}\}_{n\ge 1}$ is also an orthogonal base of  $\HH_{1}$ satisfiyng $\n{\widehat{\be}_{n}}_{\HH_{1}}^{2}=\alpha_{n}^{2}$ although $\n{\widehat{\be}_{n}}_{\HH_{0}}^{2}=1$. Morevoer the embedding $\cJ_{\hbox{\small Id}}$ is a Hilbert-Schmidt operator
\vspace*{-.1cm}
$$
\n{\cJ_{\hbox{\footnotesize Id}}}_{\cL_{2}(\HH_{0};\HH_{1})}^{2}
=\sum_{n\ge 1}\pe{\cJ_{\hbox{\footnotesize Id}}\widehat{\be}_{n}}{\cJ_{\hbox{\footnotesize Id}}\widehat{\be}_{n}}_{\rH_{1}}=\hbox{trace }\cJ_{\hbox{\footnotesize Id}}\cJ_{\hbox{\footnotesize Id}}^{*}=\sum_{n\ge 1}\alpha_{n}^{2}.
$$
So that the relative cylindrical process can be viewed as the $\cQ^{\cJ_{\hbox{\footnotesize Id}}}=\cJ_{\hbox{\footnotesize Id}}\cJ_{\hbox{\footnotesize Id}}^{*}$-Wiener process on $\HH_{0}$ taking values in a large Hilbert space $\HH_{1}$
\begin{equation}
\rW_{t}=\sum_{n\ge 1}\rB_{t}^{n}\widehat{\be}_{n},\quad t\ge 0, 
\label{eq:cylindricalWienerconstructed}
\end{equation}
(see \eqref{eq:Winerisometrygeneral}). Here 
$$
\cQ^{\cJ_{\hbox{\footnotesize Id}}}:~(\HH_{1},\n{\cdot}_{\HH_{1}})\stackrel{\cJ_{\hbox{\footnotesize Id}}^{*}}{\rightarrow}(\HH_{0},\n{\cdot}_{\HH_{0}})\stackrel{\cJ_{\hbox{\footnotesize Id}}}{\rightarrow}(\HH_{1},\n{\cdot}_{\HH_{1}})
$$
verifies $\cQ^{\cJ_{\hbox{\footnotesize Id}}}\widehat{\be}_{n}=\widehat{\be}_{n}$ and thus it is a kind of diagonal operator. Moreover
\begin{equation}
\n{\cQ^{\cJ_{\hbox{\footnotesize Id}}}}_{\cL_{1}^{+}(\HH_{1})}=\n{\cJ_{\hbox{\small Id}}}_{\cL_{2}(\HH_{0};\HH_{1})}^{2}
=\hbox{trace }\cJ_{\hbox{\footnotesize Id}}\cJ_{\hbox{\footnotesize Id}}^{*}=\sum_{n\ge 1}\alpha_{n}^{2}.
\label{eq:trazaconstructed}
\end{equation}
\fineqnum
\end{rem}
\begin{rem}\rm \label{rem:cilindro}For practical purposes, in this paper we have considered $\rH=\rL^{2}(\rI)$ and 
$$
\rW_{t}=\sum_{n\ge  1}\dfrac{1}{\sqrt{\mu_{n}}}\rB^{n}(t)\be_{n},\quad t\ge 0,
$$
(see \eqref{eq:Hypo Wiener cylindrical}) where $\{\be_{n}\}_{n\ge 1}$ and $\{\mu_{n}\}_{n\ge 1}$ was given in \eqref{eq:baseA}. Then the simplest  {\em Wiener isometry} becomes
$$
\EE\big [\n{\rW(\cdot,t)}^{2}_{\rH}\big ]=t\sum_{n\ge 1}\dfrac{1}{\mu_{n}}<+\infty.
$$
\fineqnum
\end{rem}
As  $\rG\in \cL_{2}(\HH_{1};\HH)$ implies $\rG\cJ\in \cL_{2}(\HH_{0};\HH)$, we deduce
\begin{equation}
\begin{array}{ll}
\disp \EE \big  [\big \Vert \rG\sum_{n\ge 1}\rB_{t}^{n}\cJ\widehat{\be}_{n} \big  \Vert ^{2}_{\HH}\big   ]&\hspace*{-.2cm}\disp =\lim_{n\rightarrow \infty}\EE \big  [\big \Vert \rG\sum_{m=1}^{n}\rB_{t}^{n}\cJ\widehat{\be}_{n} \big  \Vert ^{2}_{\HH}\big  ]
=\lim_{n\rightarrow \infty}\EE \big  [\big \Vert \sum_{m=1}^{n}\rB_{t}^{n}\rG\cJ\widehat{\be}_{n} \big  \Vert ^{2}_{\HH}\big  ]\\
& \disp =t\lim_{n\rightarrow \infty}\EE \big  [\big \Vert \sum_{m=1}^{n}\n{\rG\cJ\widehat{\be}_{n} }^{2}_{\HH}\big  ]=t\n{\rG\cJ}_{\cL_{2}(\HH)}^{2}\\
& \disp =t~\hbox{trace }\rG\cJ\cJ^{*}\rG^{*}<\infty.
\end{array}
\label{eq:cylindricalWienerisometry}
\end{equation}
Hence,  for any $t>0$, the series in \eqref{eq:cylindricalWienerisometry} converges in $\rL^{2}(\Omega,\cF,\PP;\HH)$ and uniformely in  $[0,\rT]~\PP-a.s.$ to a Gaussian random variable 
$$
\rG\rW_{t}=\rG\sum_{n\ge 1}\rB_{t}^{n}\cJ\widehat{\be}_{n} 
$$ 
with mean $0$ and covariation operator
$$
\cQ_{t}^{\rG\cJ}=t\hbox{ trace }\rG\cJ\cJ^{*}\rG^{*}.
$$
\par
Given $\rT>0$, we consider the set $\cP_{\rT}$ of the $\cL_{2}^{\cQ^{\cJ}}$-valued  predictable process: i.e, processes such that
$$
\rG:[0,\rT]\times \Omega \rightarrow  \cL_{2}(\HH_{1};\HH)
$$
is $\cB([0,\rT])\otimes \cF$-measurable related to $\cQ^{\cJ}$.
In fact, $\cP_{\rT}$ is a separable Hilbert space 
equipped  with  $ \n{\rG}_{\cP_{\rT}}^{2}=\cQ_{\rT}^{\rG\cJ}$ where 
\begin{equation}
\cQ_{t}^{\rG\cJ}=\int^{t}_{0}\hbox{trace }\rG_{s}\cJ\cJ^ {*}\rG_{s}^{*}ds =\int^{t}_{0}\n{\rG_{s}\cJ}_{\cL_{2}(\HH)}^{2}ds.
\label{eq:covarianceoperator}
\end{equation}
By means of density reasoning on the elementary process, as usual in integration theory, given $\rG\in\cP_{\rT}$  one introduces the stochastic integral, in the class $\cP_{\rT}$,
$$
(\rG\cdot\rW)_{t}\doteq \int^{t}_{0}\rG_{t}d \rW_{t},\quad t\in [0,\rT],
$$
satisfying
\begin{equation}
\EE \left  [\left \Vert (\rG\cdot\rW)_{t}\right \Vert ^{2}_{\HH}\right  ]=\cQ_{t}^{\rG\cJ},
\label{eq:integralisometry}
\end{equation}
(see \eqref{eq:covarianceoperator}). Then $\{(\rG\cdot \rW)_{t}\}_{t\in [0,\rT]}$ is a continuous, square integrable $\HH$-valued martingale whose covariation operator is
$$
\big [\rG\cdot \rW]_{t}=\cQ_{t}^{\rG\cJ}.
$$
In fact, \eqref{eq:integralisometry} defines an isometry 
$$
\cM^{2}_{\rT}\rightarrow \cP_{\rT},
$$
where $\cM^{2}_{\rT}$ is the space of the continuous, square integrable $\HH$-valued martingales. Moreover, the Ito formula leads to
\begin{equation}
\left \Vert \int^{t}_{0}\rG_{s}d\rW_{s}\right \Vert ^{2}_{\HH}=\int^{t}_{0}\hbox{trace }\rG_{s}\cJ\cJ^ {*}\rG_{s}^{*}ds+
\int^{t}_{0}\pe{\int^{s}_{0}\rG_{r}d\rW_{r}}{\rG_{s}d\rW_{s}}_{\HH}	ds,
\label{eq:Itoformula}
\end{equation}
which implies that $\cQ_{t}^{\rG\cJ}$ is the quadratic variation associated to the martingale $\disp t\mapsto
\int^{t}_{0}\rG_{s}d\rW_{s}$. It is easy to see that this implies the Wiener isometry, also called as Ito isometry.
\begin{rem}\rm \label{rem:extensionintegral}
\
\begin{enumerate}
	\item 	In our study of stochastic energy balance models we are choosing as coefficients suitable deterministic processes $\{\rG_{t}\}_{t\in [0,\rT]}$.
	\item With the notation of Remark \ref{rem:constructionJ} one has
	\begin{equation}\label{eq:descopositionintegral}
	\int_{0}^{t}\rG_{s}d\rW_{s}=\sum_{n\ge 1}\alpha_{n}\int_{0}^{t}\rG_{s}\widehat{\be}_{n}d\rB_{s}^{n},\quad t\in [0,\rT],
	\end{equation}
	where the series in \eqref{eq:descopositionintegral} converges  in $\rL^{2}(\Omega,\cF,\PP;\HH)$.
	
	\item One may extend the stochastic integral to the predictable process satisfying 
	$$
	\PP \left (\int^{\rT}_{0}\hbox{trace }\rG_{s}\cJ\cJ^{*}\rG_{s}^{*}ds<\infty \right )=1.
	$$
	\fineqnum
\end{enumerate}	
\end{rem}

\subsection{Some few remarks on the stochastic convolution}
\label{Sec:stochasticconvolution}

In order to simplify the exposition, here we will consider the simplified expression 
\begin{equation}
u_{t}=\rS(t)u_{0}+\rW_{t}^{\rA,\rG},\quad t>0,
\label{eq:mildsolutionshort}
\end{equation}
with respect to the identity \eqref{eq:mildsolution} (with values in $\rH=\rL^{2}(\rI)$). We want to give a precise sense, mainly, to the {\em stochastic convolution} term  
\begin{equation}
\rW_{t}^{\rA,\rG}=\int^{t}_{0}\rS(t-s)\rG_{s}d\rW_{s},\quad t\ge 0.
\label{eq:stochasticconvolution}
\end{equation}
As it was pointed out in Remark \ref{rem:semigroup}, the leading part of $(\rE_{\beta ,\rG})$ can be written as an special case of the abstract problem   
\begin{equation}
\left\{ 
\begin{array}{l}
\dfrac{d u}{d t}(t)+\rA u(t)=0,\quad t>0,\\[0.2cm] 
u(0)=u_{0}\in\rH,
\end{array}
\right.  
\label{eq:deterSellerproblemleadingf0}
\end{equation}
where the operator $\rA:~\rD(\rA)\rightarrow \rH$, introduced in \eqref{eq:leadingpartHoperator},  generates the compact semigroup $\{\rS(t)\}_{t\ge 0}$ of contractions.
\begin{rem}\rm As already said in the Introduction, the space $\rL^{2}(\rI)$ admits a Hilbertian base given by the eigenvectors of the operator $\rA$ (coinciding with the Legendre orthonormal polynomials of degree $n$)
$$
\be_{n}(x)=\sqrt{\dfrac{2}{2n+1}}\rP_{n}(x),\quad -1\le x\le 1,
$$
with $\mu_{n}=n(n+1)$ their relative eigenvalues (see \eqref{eq:baseA}). So that, for every $u_{0}\in \rL^{2}(\rI)$ the mild solution of \eqref{eq:deterSellerproblemleadingf0} admits the representation
\begin{equation}
u(t)=\rS(t)u_{0}=\sum_{n\ge 1}\pe{u_{0}}{\be_{n}}e^{-\mu_{n}t}\be_{n},\quad t\ge 0.
\label{eq:deterSellerproblemleadingsolution}
\end{equation}
\end{rem}
Next we focus on the non-autonomous problem
\begin{equation}
\left\{ 
\begin{array}{l}
\dfrac{d u}{d t}(t)+\rA u(t)=f(t),\quad t>0,\\[0.2cm] 
u(0)=u_{0}\in\HH.
\end{array}
\right. 
\label{eq:deterSellerproblemleadingf}
\end{equation}
When $u_{0}$ and $f\in\rL^{1}\big (0,\rT;\rD(\rA)\big )$ we may solve \eqref{eq:deterSellerproblemleadingf} via the generalized {\em Duhamel formula} (or constants variations formula)
$$
u(t)=\rS(t)u_{0}+\int^{t}_{0}\rS(t-s)f(s)ds,\quad 0\le s<\rT,
$$
and $u$ is called as the {\em mild solution} of \eqref{eq:deterSellerproblemleadingf}. So, by replacing \eqref{eq:deterSellerproblemleadingf} by the Ornstein-Uhlenbeck type problem
\begin{equation}
\left\{ 
\begin{array}{l}
d u_{t}+\rA u_{t}=\rG_{t}d \rW_{t},\quad t>0,\\[0.15cm] 
u(0)=u_{0}\in\HH,
\end{array}
\right. 
\label{eq:deterSellerproblemleadingRuido}
\end{equation}
the {\em mild solution} obtained by the constants variations formula leads to  \eqref{eq:mildsolutionshort}, but we must give a precise sense the so called stochastic convolution term \eqref{eq:stochasticconvolution}
in a probability space $(\Omega,\cF,\PP)$ adapted to some filtration $\{\cF_{t}\}_{t\ge 0}$.
\begin{rem}\rm
\
\begin{enumerate}
	\item For each $t>0$, we are using the notation of the process for the time-dependent terms
	$$
	\rW^{\rA,\rG}_{t}:~\Omega\rightarrow \RR. 
	$$
	When studying specifically some properties of the stochastic dynamical system the sample space $\Omega$ is chosen as an adequate space of functions (see the Section  \ref{Sec:applications}). Here, for the moment $\Omega$ is an arbitrary set.$\fin$
\end{enumerate}
\end{rem}
By simplicity, as in Remark \ref{rem:cilindro}, we consider $\rH=\rL^{2}(\rI)$ and 
$$
\rW_{t}=\sum_{n\ge  1}\dfrac{1}{\sqrt{\mu_{n}}}\rB^{n}(t)\be_{n},\quad t\ge 0,
$$
(see \eqref{eq:Hypo Wiener cylindrical}) where $\{\be_{n}\}_{n\ge 1}$ and $\{\mu_{n}\}_{n\ge 1}$ was given in \eqref{eq:baseA}. 
On the other hand, let $\rG_{s}\in \cL_{2}(\rH)$, satisfiyng \eqref{eq:infinitestochasticnoiseG} and \eqref{eq:Gpredictable} as well as 
$$	
\int_{0}^{\rT}\n{\rG_{t}\cJ} _{\cL_{2}(\rH)}^{2}dt<+\infty 	
$$
(see \eqref{eq:paraisometria}), for which the stochastic integral
$$
(\rG\cdot\rW)_{t}\doteq \int_{0}^{t}\rG_{s}d\rW_{s} 
$$
is well defined and satisfies {\em Wiener isometry} 
$$
\EE\big [\n{(\rG\cdot\rW)_{t}}^{2}_{\rL^{2}(\rI)}\big ]=\int_{0}^{t}\n{\rG_{t}\cJ} _{\cL_{2}(\rH)}^{2}dt<+\infty,
$$
(see \eqref{eq:Wienerisometryintegral}), where $\cJ:\rH\rightarrow \rH$ is given by the property $\cJ\be_{n}=\dfrac{1}{\sqrt{\mu_{n}}}\be_{n}$. 
\par
As it was proved in \cite[Theorem 0.3]{BvN}, the problem \eqref{eq:deterSellerproblemleadingRuido} has a weak solution $\{u_{t}\}_{t\in [0,\rT]}$ on $[0,\rT]$. When $u_{0}=0$ the solution is unique and it is given by the stochastic convolution
$$
u_{t}=\rW^{\rA,\rW}_{t},\quad t\in [0,\rT].
$$ 
\begin{rem}[\cite{BvN,DaPZer,DaP}]\rm \
\begin{enumerate}
	\item For any $\rT>0$ the process $\{\rW_{t}^{\rA,\rG}\}_{t\in [0,\rT]}$ is a mean square continuos adapted process taking values on $\HH$. Thus it belongs to the Banach space
	$$
	\cC_{\rW}\big (\big [0,\rT\big ];\rL^{2}(\Omega,\cF,\PP;\rH)\big ) \hbox { or  $ \cC_{\rW}\big (\big [0,\rT\big ];\rH)$, in short},
	$$
	of the continuous processes $\{\rX_{t}\}_{t\in [0,\rT]}$ adapted to  $\{\rW_{t}\}_{t\in[0,\rT]}$ such that $\rX_{t}$ is measurable
	endowed with the norm
	$$
	\n{\rX}_{\cC_{\rW}([0,\rT];\rH)}\doteq \left (\sup_{t\in [0,\rT]}\EE\big [\n{\rX_{t}}^{2}_{\rH}\big ]\right )^{\frac{1}{2}}.
	$$
	\item The process $\{\rW_{t}^{\rA,\rG}\}_{t\in [0,\rT]}$ is $\PP$-almost surely continuous on $[0,\rT]$.
	\item The {\em Wiener isometry }
	\begin{equation}
	\EE   \big [\n{\rW_{t}^{\rA,\rG}}^{2}_{\rH}\big ]=\cQ_{t}^{\cA,\rG}<\infty,\quad t\ge 0,
	\label{eq:integralisometryconvolution}
	\end{equation}
	holds for the covariation operator
	\begin{equation}
	\cQ_{t}^{\rA,\rG}=\int^{t}_{0}\hbox{trace }\rS(t-s)\rG\cQ^{\cJ_{\hbox{\footnotesize Id}}}\rG^{*}\rS(t-s)^{*}ds.
	\label{eq:covarianceoperatorconvolution}
	\end{equation}
	\fineqnum
\end{enumerate}
\end{rem}
\begin{rem}\rm In the case of Legendre diffusion we have the property
$$
\int^{t}_{0}\hbox{trace }\rS(s)\cQ^{\cJ_{\hbox{\small Id}}}\rS(s)^{*}=
\int^{t}_{0}\sum_{n\ge 1}\dfrac{e^{-2s\mu_{n}}}{\mu_{n}}ds=\dfrac{1}{2}\sum_{n\ge 1}\dfrac{1-e^{-2t\mu_{n}}}{\mu_{n}^{2}}\le \dfrac{1}{2}\sum_{n\ge 1}\dfrac{1}{\mu_{n}^{2}}<+\infty,
$$
since $\mu_{n}=n(n+1).\fin$
\end{rem}

Recently, a new {\em Burkholder-Rosental type inequality} 
\begin{equation}
\EE \left [\sup_{t\in[0,\rT]}\n{\rW_{t}^{\rA,\rG}}_{\rH}^{p}\right ]\le \rC_{p}\left (\int^{t}_{0}\hbox{trace }\rS(t-s)\rG\cQ^{\cJ_{\hbox{\footnotesize Id}}}\rG^{*}\rS(t-s)^{*}ds\right )^{\frac{p}{2}},\quad 2\le p<\infty,
\label{eq:Burkholder-Rosental}
\end{equation}
was obtained in \cite{NJvV}, where the positive constant $\rC_{p}$ is the order $O(\sqrt{p})$. It coincides with the {\em Burkholder inequality} (and therefore it is optimal) as $p\rightarrow \infty$. Among other properties, the inequality \eqref{eq:Burkholder-Rosental} enables us to obtain the stablity and pointwise uniform convergence of some time discretisation schemes (see \cite{NJvV} for details).

\bigskip

\begin{tabular}{ll}
Gregorio Díaz & Jesús Ildefonso  Díaz \\
& Instituto Matemático Interdisciplinar (IMI)\\
Dpto. Análisis Matemático & Dpto. Análisis Matemático \\
y Matem\'atica Aplicada & y Matem\'atica Aplicada \\
U. Complutense de Madrid & U. Complutense de Madrid \\
28040 Madrid, Spain & 28040 Madrid, Spain\\
{\tt gdiaz@ucm.es} & {\tt jidiaz@ucm.es}
\end{tabular}

\end{document}